\documentclass[preprint,12pt]{article}
\usepackage{amssymb}
\usepackage{mathrsfs}
\usepackage{amsfonts}
\usepackage{graphicx}
\usepackage{colortbl,dcolumn}
\usepackage{amsmath}
\usepackage{psfrag}
\usepackage{booktabs}
\usepackage{slashbox}
\usepackage{cite}
\numberwithin{equation}{section}
\usepackage[top=1.2in, bottom=1.2in, left=1in, right=1in]{geometry}

\newtheorem{theorem}{Theorem}[section]
\newtheorem{definition}{Definition}
\newtheorem{lemma}{Lemma}

\newtheorem{remark}{Remark}
\newtheorem{example}{Example}

\newtheorem{scheme}{\textbf{Scheme}}
\begin{document}
\title{Near Preservation of Quadratic Invariants by \\Stochastic Runge-Kutta Methods \footnotemark[1]}
       \author{
        Jialin Hong\footnotemark[2], Lijin Wang\footnotemark[3], Dongsheng Xu\footnotemark[4] and Liying Zhang\footnotemark[1]\\
       {\small \footnotemark[2]~\footnotemark[4]~\footnotemark[1] Institute of Computational Mathematics and Scientific/Engineering Computing,}\\{\small Academy of Mathematics and Systems Science, Chinese Academy of Sciences, }\\
         {\small Beijing 100190, P.R.China } \\
       {\small \footnotemark[3] School of Mathematical Sciences, University of Chinese Academy of Sciences,}\\
    {\small Beijing 100190, P.R.China }}
       \maketitle
       \footnotetext{\footnotemark[2] The firs authors is supported by the NNSFC (No. 91130003, No. 11021101 and No. 11290142).\\
        \footnotemark[3]The second author is supported by the NNSFC (No. 11071251). \\\footnotemark[1]The fourth author is supported by the NNSFC (NO. 11301001, NO. 2013SQRL030ZD).}
        \footnotetext{\footnotemark[1]Corresponding author: lyzhang@lsec.cc.ac.cn}

       \begin{abstract}
          {\rm\small Based on the combinatory theory of rooted colored trees, we investigate the conditions for the explicit stochastic Runge-Kutta (SRK) methods to preserve quadratic invariants (QI) up to certain orders of accuracy.
These conditions can supply a practical approach of constructing explicit nearly conservative SRK methods. Meanwhile, we estimate errors in the preservation of QI resulting from iterative implementation of implicit conservative SRK methods with fixed-point and Newton's iterations.
Finally, numerical experiments are performed to test the behavior of the methods in preserving QI.}\\


\end{abstract}

\section{Introduction}
\label{sect:intro}

\label{sec;introduction}
Stochastic differential equations (SDEs)
can describe the natural and social phenomena more realistically than deterministic differential equations. However, they are usually very difficult to be solved analytically, which gives rise to the research on numerical methods for SDEs (see \cite{Higham,2,3} and references therein).

Designing numerical methods inheriting qualitative properties of the original SDE systems is an attractive topic of research. As one of the important qualitative properties, the invariants of the underlying continuous differential equations systems are expected to be preserved by numerical methods for the reliability and a long time stability.
There are many references concerning this problem. \cite{Schurz} studies the midpoint (trapezoidal) methods preserving the first and second moments for linear SDEs; \cite{Misama} proposes the conserving energy difference scheme for stochastic dynamical systems; \cite{Milstein} constructs symplectic numerical schemes for stochastic canonical Hamiltonian problems;  \cite{8,9} develop the generating functions for stochastic symplectic methods; \cite{7} investigates the boundary preserving semianalytic numerical algorithms for SDEs; \cite{15} designs invariants-preserving methods for SDEs by using the discrete gradient approach;
The recent works
\cite{1} and \cite{Chen} propose a new class of energy-preserving numerical schemes for stochastic Hamiltonian systems with non-canonical structure matrix and present a novel conservative method for numerical computation of
general stochastic differential equations in the Stratonovich sense with a conserved quantity, respectively.

Quadratic invariants (QI) cover a large class of important properties, such as the symplecticity of the original systems. Stochastic Runge-Kutta (SRK) methods preserving QI, the so called conservative SRK methods, are studied in \cite{15}. It is proved that all SRK methods that preserve QI must satisfy certain conditions which indicate that they are in general fully implicit. Therefore these implicit methods build up a barrier for the implementation. On the other hand, explicit SRK methods are easily realized, while fail to preserve QI accurately. Currently, there are no references concerning the problem of how to balance the computation complexity and the preservation
of QI, that is, to find certain methods that are explicit but can preserve QI up to satisfying order.
So based on the colored rooted tree theory, we give in this paper the conditions that guarantee the preservation of QI up to any desired order of accuracy.
 Meanwhile, using these conditions we construct two explicit schemes facilitating the numerical implementation in the last parts.

Another alternative of removing implicitness of the conservative SRK methods is to perform iterations on the implicit methods.
A direct consequence of such treatments is the loss of the QI-preservation. In this paper, we estimate how far the QI-preservation will be ruined by the iterations by giving the error bounds which is related to both iteration numbers and time step-sizes. This gives hints to appropriate choices of the parameters of methods.

The paper is organized as follows. Section 2 investigates the near preservation of quadratic invariants by explicit SRK methods, while that by iterative implementation of implicit SRK methods is studied in Section 3. Section 4 performs numerical experiments, succeeded by concluding remarks in Section 5.

In the sequel, we will make use of the following notation:
$\|x\|$ is the Euclidean norm of a vector x or the induced norm for a matrix.
\section{Near Preservation of Quadratic Invariants by Explicit SRK Methods}
\label{sec:1}
Let $(\Omega, \mathcal{F}, \{\mathcal{F}_{t}\}_{t \geq 0}, \mathbf{P})$ be a complete probability space with a filtration $\{\mathcal{F}_{t}\}_{t\geq 0}$,
which is a nondecreasing right continuous family of $\sigma$ subalgebra of $\mathcal{F}$, and where $\mathcal{F}_0$ contains all the $\mathbf{P}$-null sets in $\mathcal{F}$.
 Let $W(t)$ be a $\mathcal{F}_{t}$ adapted one-dimensional standard Brownian motion.

Consider the $d$-dimensional autonomous SDE of Stratonovich sense
\begin{equation}\label{SDE}
dy=f(y)dt+g(y)\circ dW(t),\quad y(0)=y_0,~~ t\in[0,T],
\end{equation}
where $y_{0}$ is independent of the $\sigma$-algebra generated by the Brownian motion, $E\|y_{0}\|^{2}<\infty$.
In fact, based on the relation between the It\^o and Stratonovich integrals, the solution of (\ref{SDE}) is also the solution of the It\^o SDE
\begin{equation}\nonumber
dy=\underline{f}(y)dt+g(y)dW(t),\quad y(0)=y_0,~~ t\in[0,T],
\end{equation}
where $\underline{f}(y)=f(y)+\frac 1 2\frac{\partial g(y)}{\partial y}g(y)$. Additionally, under the assumption that the drift $\underline{f}:  R^d\rightarrow R^d$ and the diffusion
$g: R^d\rightarrow R^d$ are measurable functions satisfying
 \begin{equation}\nonumber
\|\underline{f}(y)-\underline{f}(x)\|+\|g(y)-g(x)\|\leq L_1\|y-x\|,
\end{equation}
\begin{equation}\nonumber
\|\underline{f}(y)\|^2+\|g(y)\|^2\leq L_2(1+\|y\|^2),
\end{equation}
where $x, y\in R^d$ and $L_1, L_2$ are positive constants, the solution of (\ref{SDE}) exists and is unique according to \cite{90}.
\begin{definition}\cite{4}
A differentiable scalar function $I(y)$ is called an invariant of (\ref{SDE}) if it satisfies
\begin{equation}\label{def}
\nabla I^{T}f(y)=0,~~\nabla I^{T}g(y)=0,
\end{equation}
where $\nabla I=(\frac{\partial I}{\partial y^{1}},\cdots,\frac{\partial I}{\partial y^{d}})^{T}$.

Specially, if  $I(y)=y^{T}Cy$, where $C$ is a symmetric square matrix, then $I(y)$ is said to be a quadratic invariant of (\ref{SDE}).
\end{definition}

We consider the following class of stochastic Runge-Kutta methods as in \cite{10}
\begin{eqnarray}\label{mysrk}
Y_i&=&y_n+h\sum_{j=1}^s a_{ij}f(Y_j)+\Delta W_{n}\sum_{j=1}^sb_{ij}g(Y_{j}),\quad i=1,\cdots,s,\nonumber\\
y_{n+1}&=&y_{n}+h\sum_{i=1}^s \alpha_i f(Y_i)+\Delta W_n\sum_{i=1}^s \beta_i g(Y_i),\quad n=0,\cdots, N_T-1.
\end{eqnarray}
 Introduce the notations
 \begin{eqnarray*}
 &e^T=(\underset{s}{\underbrace{1,\cdots,1}}),\quad A=(a_{ij}),\quad B=(b_{ij}), \\
 &\alpha^T=(\alpha_1,\cdots,\alpha_{s}), \quad \beta^T=(\beta_{1},\cdots,\beta_{s}),\quad Y=(Y_{1}^{T},\ldots,Y_{s}^{T})^{T}, \\
& F(Y)=(f(Y_{1})^{T}, \ldots,f(Y_{s})^{T})^{T}, \quad G(Y)=(g(Y_{1})^{T},\ldots,g(Y_{s})^{T})^{T}.
 \end{eqnarray*}
 Then the stochastic Runge-Kutta method (\ref{mysrk}) can be rewritten as
\begin{eqnarray}\label{SRK-1}
Y&=&e\otimes y_{n}+h(A\otimes I_d)F(Y)+\Delta W_n(B\otimes I_d)G(Y),\nonumber\\
y_{n+1}&=&y_{n}+h(\alpha^{T}\otimes I_d)F(Y)+\Delta W_n(\beta^{T}\otimes I_d)G(Y),
\end{eqnarray}
for $n=0,\cdots,N_T-1$, where $\otimes$ denotes the Kronecker (tensor) product. The Butcher tableau for (\ref{SRK-1}) is
\begin{equation}
\begin{tabular}{|cc}
 $A$& ~~$B$ \\
  \hline
   $\alpha^T$ & ~~$\beta^T$  \\
\end{tabular}.\nonumber
\end{equation}

\begin{theorem}\label{thm1}(\cite{15}) If the coefficients of the SRK method (\ref{SRK-1}) satisfy
\begin{eqnarray}\label{conditionnew}
\alpha_{i}a_{ij}+\alpha_{j}a_{ji}&=\alpha_{i}\alpha_{j},\\\nonumber
\alpha_{i}b_{ij}+\beta_{j}a_{ji}&=\alpha_{i}\beta_{j},\\\nonumber
\beta_{i}b_{ij}+\beta_{j}b_{ji}&=\beta_{i}\beta_{j},
\end{eqnarray}
for all $i,j=1,\cdots,s$, then it preserves quadratic invariants of (\ref{SDE}).
\end{theorem}

By convention, SRK methods preserving quadratic invariants are called conservative SRK methods.

\begin{remark} \label{rm1}The $2d$-dimensional stochastic Hamiltonian system
\begin{eqnarray}\label{SHS}
d p(t)&=&-\frac{\partial H(p,q)}{\partial q}dt-\frac{\partial H^1(p,q)}{\partial q}\circ d W(t),~~p(0)=p_0,\nonumber \\
d q(t)&=&~~\frac{\partial H(p,q)}{\partial p}dt+\frac{\partial H^1(p,q)}{\partial p}\circ d W(t),~~q(0)=q_0,
\end{eqnarray}
with $t\in[0,T]$, where $p$, $q\in R^d$, $H$, $H^1$ are differentiable scalar functions, which also takes the equivalent form
\begin{equation}\label{SHE-equivalent}
dy=J^{-1}\nabla H(y)dt+J^{-1}\nabla H^1(y)\circ dW(t),\quad y(0)=y_0,
\end{equation}
with $y=(p^T,q^T)^T$ and \[J=
   \left(\begin{array}{ccc}
   0&I_d\\
   -I_d&0
   \end{array}\right).
   \]
The equations (\ref{SHE-equivalent}) possesses the invariant  $\psi^TJ\psi$, with $\psi=\frac{\partial y}{\partial y_0}$, i.e.,
\begin{equation}\label{symplectic-2}
\left(\frac{\partial y(t)}{\partial y_{0}}\right)^{T}J\left(\frac{\partial y(t)}{\partial y_{0}}\right)=J,\quad \forall\,\, 0\leq t\leq T,
\end{equation}
which is equivalent to
\begin{equation}\label{symplectic-1}
dp(t)\wedge dq(t)=d p_0 \wedge d q_0,\quad \forall \,\,0\leq t\leq T,
\end{equation}
namely, the symplecticity. Numerical methods that preserve the symplecticity are called symplectic methods, with the characterization
\begin{equation}
dp_{n+1}\wedge dq_{n+1}=dp_{n}\wedge dq_{n},\quad \forall\,\,0\leq n\leq N_T-1,
\end{equation}
or equivalently
\begin{equation}
dy_{n+1}\wedge J dy_{n+1}=dy_{n}\wedge J dy_{n},\quad \forall\,\,0\leq n\leq N_T-1.
\end{equation}

It is proved in \cite{15} that, the symplectic structure $\psi^T
J\psi$ is a quadratic invariant of the following augmented system of (\ref{SHE-equivalent})
\begin{align}\label{aug}
dy&=J^{-1}\nabla H(y)dt+J^{-1}\nabla H^1(y)\circ dW(t),\quad y(0)=y_0,\\
d\psi&=J^{-1}\nabla^2 H(y)\psi dt+J^{-1}\nabla^2 H^1(y)\psi\circ
dW(t),\quad \psi(0)=I_d.
\end{align}
Thus, as a consequence of Theorem \ref{thm1}, each conservative SRK method applied to the stochastic Hamiltonian system (\ref{SHE-equivalent}) is a symplectic SRK method, as proved in \cite{15}.
\end{remark}

From the conditions (\ref{conditionnew}) we can see that, all SRK methods preserving quadratic invariants are in general fully implicit. It is worthy to consider explicit SRK methods that preserve quadratic invariants up to a certain order when applying to practical implementation. A numerical discretisation $\{y_k\}_{k=0}^{N_T}$ is said to have order $\flat$ of quadratic invariants conservation if system (\ref{SDE}) possessing the quadratic invariant $I(y)=y^{T}Cy$ holds $I(y_{N_T})-I(y_0)=O(h^{\flat})$ for all $y_0$ a.s.

We start by recalling some notations and properties of rooted colored trees (see \cite{10}). Since integration in stochastic case is with respect to $dt$ and $dW(t)$,  each node of a tree can be colored with any one of the two colorings $\{0,1\}$. A node colored with the label $0$ corresponds to integration with respect to $dt$ and the node is called a deterministic node. A node colored with label $1$ is called a stochastic node. Let $\tau_{k} (k=0,1)$ denote the tree with a single node with color $k$ and $\Gamma$ be the set of all rooted trees with all possible colorings. let $\iota=[\iota_1,\ldots,\iota_{\ell}]_{k}$ be the tree formed by joining subtrees $\iota_1,\ldots,\iota_{\ell}$ each by a single branch to a common root with color $k$. Then the elementary differential associated with $\iota=[\iota_1,\ldots,\iota_{\ell}]_{k}$ is
\begin{displaymath}
F(\iota)y=g_{k}^{(\ell)}(F(\iota_1)y,\ldots,F(\iota_{\ell})y),~~~F(\phi)y=y.
\end{displaymath}

Let $\circ$ denote the deterministic node and  $\bullet$ the stochastic node, respectively. Let $\Gamma_0$ be the collection of all trees with deterministic root $\circ$, and $\Gamma_1$ be that with stochastic root $\bullet$, respectively.

\begin{theorem}\label{mainthm1}
Suppose that a $s-$stage stochastic Runge-Kutta methods (\ref{SRK-1}) is applied to the stochastic system (\ref{SDE}) with a quadratic invariant $I(y)=y^{T}Cy$, then
\begin{eqnarray}
I(y_{n+1})-I(y_{n})&=&-\sum_{\iota\in \Gamma_0, \iota'\in \Gamma_{0}}\omega(\iota)\omega(\iota')(h\Phi(\iota))^{T}M^{0}(h\Phi(\iota))\Omega(\iota,\iota')(y_n)\nonumber\\
&-&\sum_{\iota\in \Gamma_0, \iota'\in \Gamma_{1}}\omega(\iota)\omega(\iota')(h\Phi(\iota))^{T}M^{*}(\Delta W_n \Phi(\iota))\Omega(\iota,\iota')(y_n)\nonumber\\
&-&\sum_{\iota'\in \Gamma_1,\iota'\in \Gamma_{1}}\omega(\iota)\omega(\iota')(\Delta W_n\Phi(\iota))^{T}M^{1}(\Delta W_n \Phi(\iota))\Omega(\iota,\iota')(y_n),\nonumber
\end{eqnarray}
for $n=0,\cdots,N_T-1$, where $\omega(\iota)=\alpha(\iota)\gamma(\iota)/\rho(\iota)!$ and $\Omega$ is a bilinear form on the elementary differentials defined by $\Omega(\iota,\iota')(y_n)=F(\iota)(y_n)^{T}CF(\iota)(y_n)$, and
\begin{eqnarray}
M^0=diag(\alpha_1,\ldots,\alpha_s)A+A^{T}diag(\alpha_1,\ldots,\alpha_s)-\alpha \alpha^{T},\nonumber\\ M^{1}=diag(\beta_1,\ldots,\beta_s)B+B^{T}diag(\beta_1,\ldots,\beta_s)-\beta \beta^{T},\nonumber\\   M^{*}=diag(\alpha_1,\ldots,\alpha_s)B+A^{T}diag(\beta_1,\ldots,\beta_s)-\alpha \beta^{T}.\nonumber
\end{eqnarray}
\end{theorem}
$Proof.$
Set $f_i=f(Y_i)$ and $g_i=g(Y_i)$. From $(\ref{SRK-1})$, we have
\begin{eqnarray}
I(y_{n+1})-I(y_n)&=&-\sum_{i,j=1}^{s}m_{ij}^{0}(h f_i)^{T}C(h f_j)-2\sum_{i,j=1}^{s}m_{ij}^{*}(h f_i)^{T}C(\Delta W_n g_j)\nonumber \\
&-&\sum_{i,j=1}^{s}m_{ij}^{1}(\Delta W_n g_i)^{T}C(\Delta W_n g_j).
\end{eqnarray}
Firstly, let
$z^{(0)}=h e_{i}$ and $z^{(1)}=0$ with $e_{i}=(0,\ldots,1,\ldots,0)^{T}$, then we have
\begin{eqnarray}
y_{n+1}&=&y_{n}+h f_i,\nonumber \\
y_{n+1}&=&y_{n}+\sum_{\iota\in \Gamma_{0}}\frac{\gamma(\iota)}{\rho(\iota)!}\alpha(\iota)h e_{i}^{T}\Phi(\iota)F(\iota)(y_n).\nonumber
\end{eqnarray}
Therefore
\begin{displaymath}
h f_i=\sum_{\iota\in \Gamma_{0}}\frac{\gamma(\iota)}{\rho(\iota)!}\alpha(\iota)h e_{i}^{T}\Phi(\iota)F(\iota)(y_n).
\end{displaymath}
Similarly, let $z^{(0)}=0$ and $z^{(1)}=\Delta W_n e_i$ with $e_{i}=(0,\ldots,1,\ldots,0)^{T}$, we obtain
\begin{displaymath}
\Delta W_n g_i=\sum_{\iota\in \Gamma_{1}}\frac{\gamma(\iota)}{\rho(\iota)!}\alpha(\iota)\Delta W_n e_{i}^{T}\Phi(\iota)F(\iota)(y_n).
\end{displaymath}
Thus we can deduce that
\begin{eqnarray}
-\sum_{i,j=1}^{s}m_{ij}^{0}(h f_i)^{T}C(h f_j)=-\sum_{\iota\in \Gamma_0, \iota'\in \Gamma_{0}}\omega(\iota)\omega(\iota')(h\Phi(\iota))^{T}M^{0}(h\Phi(\iota))\Omega(\iota,\iota')(y_n).\nonumber
\end{eqnarray}
Similarly, it holds
\begin{eqnarray*}
&-\sum_{i,j=1}^{s}m_{ij}^{*}(h f_i)^{T}C(\Delta W_n g_j)\\&=-\sum_{\iota\in \Gamma_0, \iota'\in \Gamma_{1}}\omega(\iota)\omega(\iota')(h\Phi(\iota))^{T}M^{*}(\Delta W_n \Phi(\iota))\Omega(\iota,\iota')(y_n),\nonumber
\end{eqnarray*}
and
\begin{eqnarray*}
&-\sum_{i,j=1}^{s}m_{ij}^{1}(\Delta W_n g_i)^{T}C(\Delta W_n g_j)\\&=-\sum_{\iota'\in \Gamma_1, \iota'\in \Gamma_{1}}\omega(\iota)\omega(\iota')(\Delta W_n\Phi(\iota))^{T}M^{1}(\Delta W_n \Phi(\iota))\Omega(\Gamma,\iota')(y_n).\nonumber
\end{eqnarray*}
Then the result follows immediately.

\begin{remark}\label{rm2}
For implicit SRK methods satisfying the condition (\ref{conditionnew}), $M^0$, $M^1$ and $M^{*}$ are null matrices, which imply that they preserve the quadratic invariant $I(y)$ accurately, coinciding with the result of Theorem \ref{thm1}.
\end{remark}

The following Lemma is a consequence of the Borel-Cantelli Lemma and provides a relation between the convergence rate in the $p$th mean and the path-wise convergence rate.
\begin{lemma}\label{lemma1}(\cite{16})
Let $\alpha >0$ and $K(p)\in [0,\infty)$ for $p \geq 1$. In addition, let $Z_n$, $n\in N$, be a sequence of random variables such that
\begin{displaymath}
(E|Z_n|^{p})^{1/p}\leq K(p)\cdot n^{-\alpha}
\end{displaymath}
for all $p\geq 1$ and all $n\in N$. Then for all $\epsilon >0$ there exists a finite and non-negative random variable $\eta_{\epsilon}$ such that
\begin{displaymath}
|Z_n|\leq \eta_{\epsilon}\cdot n^{-\alpha+\epsilon}~~~a.s.
\end{displaymath}
for all $n\in N$.
\end{lemma}

Given a rooted colored tree $\iota$ with $n_0$ deterministic nodes and $n_{1}$ stochastic nodes, then the order of the tree, $ord(\iota)$, is $ord(\iota)=n_0+\frac{1}{2}n_1$ (see \cite{10}). We have the following result.

\begin{theorem}\label{mainthm2}
If $\Phi(\iota)^{T}M^{0}\Phi(\iota')=0$ for all $\iota, \iota' \in \Gamma_0$ such that  $ord(\iota)+ord(\iota')\leq \gamma+1/2$, $\Phi(\iota)^{T}M^{*}\Phi(\iota')=0$ for all $\iota\in \Gamma_0, \iota' \in \Gamma_1$ such that $ord(\iota)+ord(\iota')\leq \gamma+1/2$ , and $\Phi(\iota)^{T}M^{1}\Phi(\iota')=0$ for all $\iota, \iota' \in \Gamma_1$ such that $ord(\iota)+ord(\iota')\leq \gamma+1/2$ , then for all $\epsilon >0$ and $N_Th\leq1$, there exists a finite non-negative random variable $C(\epsilon)$ such that
\begin{equation}\label{main thm2 eq}
|I(y_{N_T})-I(y_0)|\leq C(\epsilon)h^{\gamma-\epsilon},\quad a.s.,
\end{equation}
where we call $\gamma-\epsilon$ the order of preservation of quadratic invariants.
\end{theorem}
$Proof.$
From the definition of $\Phi(\iota)$ and the fact that $(E|\Delta W_n|^p)^{1/p}\leq Ch^{1/2}$, we have
\begin{displaymath}
(E|h\Phi(\iota)|^{p})^{\frac{1}{p}}\leq Ch^{ord(\iota)}
\end{displaymath}
for all $\iota\in \Gamma_{0}$,
and
\begin{displaymath}
(E|\Delta W_n\Phi(\iota)|^{p})^{\frac{1}{p}}\leq Ch^{ord(\iota)},
\end{displaymath}
for all $\iota\in\Gamma_{1}$.

By Theorem \ref{mainthm1} together with  the Minkowski inequality, we obtain
\begin{displaymath}
(E|I(y_{n+1})-I(y_{n})|^{p})^{1/p}\leq Ch^{\gamma+1},
\end{displaymath}
which, according to the Lemma \ref{lemma1}, implies that for all $\epsilon >0$, there exists a finite non-negative random variable $C(\epsilon)$ such that
\begin{displaymath}
|I(y_{n+1})-I(y_{n})|\leq C(\epsilon)h^{\gamma+1-\epsilon},~~~a.s.
\end{displaymath}
for $n=0,\cdots,N_T-1$.
This completes the proof.

\begin{remark}\label{rm3}\hfill
\begin{itemize}
\item For implicit SRK methods with (\ref{conditionnew}), the conditions of the Theorem \ref{mainthm2} are naturally satisfied, and the left-hand-side of equation (\ref{main thm2 eq}) equals zero.
\item Analogous to the definition of pseudo-symplecticity in deterministic cases, we can define stochastic pseudo-symplecticity. A one-step method $\Psi_{h}(y_n)$ applied to a stochastic Hamiltonian system with step size $h$ is called stochastic pseudo-symplectic order $\gamma$ if it holds
\begin{displaymath}
|\Psi_{h}'(y_n)^{T}J\Psi_{h}'(y_n)-J|\leq C(\epsilon) h^{\gamma-\epsilon},~~a.s.
\end{displaymath}
for $n=1,\cdots,N_T$, where $\Psi_{h}'(y_n)=\frac{\partial \Psi_{h}(y_n)}{\partial y_n}$, and $C(\epsilon)$ is a finite non-negative random variable. Then according to the discussion in Remark \ref{rm1}, each SRK method (\ref{SRK-1}) with order $\gamma-\epsilon$ of preservation of quadratic invariants is stochastic pseudo-symplectic of order $\gamma$.
\item One can expect to construct explicit SRK methods up to a certain order $\gamma-\epsilon$ by letting the coefficients of the explicit SRK methods satisfy the conditions of the Theorem \ref{mainthm2}.
\end{itemize}
\end{remark}

\begin{example}\label{e1} A SRK scheme of order $(1.0,2.0-\epsilon)$

For notational simplicity, we call the scheme is of order $(p,q-\epsilon)$ if it is of global convergence order $p$ and has order $q-\epsilon$ of preservation of quadratic invariants. In \cite{10}, general order conditions for general SRK methods (\ref{SRK-1}) have been given. A SRK method (\ref{SRK-1}) has strong global order $1.0$, if it satisfies
\begin{equation}\label{order-con}
\alpha^Te=1,~~\beta^Te=1,~~\beta^TBe=\frac{1}{2}.
\end{equation}
From the rooted colored tree theory, we have that up to order $2.5$
\begin{eqnarray*}
\Gamma_{0}=&\{\tau_{0},[\tau_{1}]_{0},[\tau_1,\tau_1]_0,[[\tau_1]_1]_0,[\tau_0]_0,[[\tau_1]_0]_0,[[\tau_0]_1]_0,[\tau_0,\tau_1]_0,\\&[\tau_1,\tau_1,\tau_1]_0,
[\tau_1,[\tau_1]_1]_0,[[[\tau_1]_1]_1]_0,[[\tau_1,\tau_1]_1]_0,\ldots\}
\end{eqnarray*}
and
\begin{eqnarray*}
\Gamma_{1}=&\{\tau_1,[\tau_1]_1,[\tau_0]_1,[[\tau_1]_1]_1,[\tau_1,\tau_1]_1,[\tau_1,[\tau_1]_1]_1,[[[\tau_1]_1]_1]_1,[[\tau_1,\tau_1]_1]_1,[\tau_1,\tau_1,\tau_1]_1,\\
& [[\tau_0]_1]_1,[\tau_1,\tau_0]_1,[[\tau_1]_0]_1,[\tau_0,\tau_0]_1,[[\tau_0]_0]_1,[[[\tau_0]_1]_1]_1,[[\tau_1,\tau_0]_1]_1,\\
&[\tau_0,[\tau_1]_1]_1,[\tau_0,\tau_1,\tau_1]_1,[\tau_1,[\tau_0]_1]_1,[[\tau_1]_1,[\tau_1]_1]_1,[\tau_1,[[\tau_1]_1]_1]_1,\\
&[\tau_1,\tau_1,[\tau_1]_1]_1, [\tau_1,[\tau_1,\tau_1]_1]_1,[[[[\tau_1]_1]_1]_1]_1,[[[\tau_1,\tau_1]_1]_1]_1,[[[\tau_1]_1,\tau_1]_1]_1,\\
&[[\tau_1,\tau_1,\tau_1]_1]_1, [\tau_1,\tau_1,\tau_1,\tau_1]_1,[[[\tau_1]_0]_1]_1,[\tau_1,[\tau_1]_0]_1,[[\tau_1,\tau_1]_0]_1,\\
&[[[\tau_1]_1]_0]_1, \ldots\}
\end{eqnarray*}

According to Theorem \ref{mainthm2}, a SRK method of order $(1.0,2.0-\epsilon)$, should satisfy, additional to (\ref{order-con}), the following conditions
\begin{eqnarray*}
e^{T}M^{*}e=0,~~e^{T}M^{1}e=0,~~e^{T}M^{1}Be=0,\nonumber\\
e^{T}M^{0}e=0,~~e^{T}M^{*}Be=0,~~(Be)^{T}M^{*}e=0,~~e^{T}M^{1}Ae=0\nonumber\\
e^{T}M^{1}(Be)^2=0,~~(Be)^{T}M^{1}(Be)=0,~~e^{T}M^{1}B^2e=0,\nonumber\\
\end{eqnarray*}
and
\begin{eqnarray*}
&e^{T}M^{0}Be=0,~~e^{T}M^{*}Ae=0,~~~e^{T}M^{*}B^2e=0,~~~e^{T}M^{*}(Be)^2=0,\\
&(Be)^{T}M^{*}Be=0,~~((Be)^{2})^{T}M^{*}e=0,~~(B^2e)^{T}M^{*}e=0,~~(A e)^{T}M^{*}e=0,\\
&e^{T}M^{1}(Be\cdot B^2e)=0,~~e^{T}M^{1}(B^3e)=0,~~e^{T}M^{1}B(Be)^2=0,\\
&e^{T}M^{1}(Be)^3=0,~~e^{T}M^{1}(BAe)=0,~~e^{T}M^{1}(Be\cdot A e)=0,\\
&e^{T}M^{1}(ABe)=0,~~(Be)^{T}M^{1}(Ae)=0,\\
&(Be)^{T}M^{1}(B^2e)=0,~~(Be)^{T}M^{1}(Be)^2=0,
\end{eqnarray*}
from which we can establish the following SRK scheme of order $(1.0,2.0-\epsilon)$.
\begin{scheme}\label{scheme1}
\begin{equation}
\begin{tabular}{|ccccccccccccccccc}
$0$& $0$ & $0$& ~~$0$& $0$ & $0$  \\
$\frac{1}{4}$&  $0$ & $0$&  ~~$\frac{1}{4}$&  $0$ & $0$  \\
$-\frac{1}{2}$&  $\frac{3}{2}$ &$0$ &~~$-\frac{1}{2}$&  $\frac{3}{2}$ &$0$  \\
\hline
$0$& $\frac{2}{3}$ & $\frac{1}{3}$&~~$0$& $\frac{2}{3}$ & $\frac{1}{3}$ \\
\end{tabular}.
\end{equation}
\end{scheme}
\end{example}

\begin{example} \label{e2}A SRK scheme of order $(1.0,2.5-\epsilon)$

Similar to the discussion above, for a SRK method of order $(1.0,2.5-\epsilon)$, the following additional conditions should be satisfied
\begin{eqnarray*}
&&e^{T}M^{0}(Be)^2=0,~~e^{T}M^{0}(B^2e)=0,~~e^{T}M^{0}(Ae)=0,~~(Be)^{T}M^{0}(Be)=0,\nonumber\\
&&e^{T}M^{*}(Be\cdot B^2e)=0,~~e^{T}M^{*}(B^3e)=0,~~e^{T}M^{*}B(Be)^2=0,\nonumber\\
&&e^{T}M^{*}(Be)^3=0,~~e^{T}M^{*}(BAe)=0,~~e^{T}M^{*}(Be\cdot Ae)=0,\nonumber\\
&&e^{T}M^{*}(ABe)=0,~~(Be)^{T}M^{*}(B^2e)=0,~~(Be)^{T}M^{*}(Ae)=0,\nonumber\\
&&(Be)^{T}M^{*}(Be)^2=0,~~(Be)^{2T}M^{*}(Be)=0,~~(B^{2}e)^TM^{*}(Be)=0,\nonumber\\
&&(A e)^{T}M^{*}(Be)=0,~~(A Be)^{T}M^{*}e=0,~~(B A e)^{T}M^{*}e=0,\nonumber\\
&&(A e\cdot B e)^{T}M^{*}e=0,~~(Be)^{3T}M^{*}e=0,~~(B e\cdot B^2 e)^{T}M^{*}e=0,\nonumber\\
&&(B^3 e)^{T}M^{*}e=0,~~(B(Be)^2)^{T}M^{*}e=0,~~e^{T}M^{1}(Ae)^2=0,\nonumber\\
&&e^{T}M^{1}(A^2e)=0,~~e^{T}M^{1}(B^2Ae)=0,~~e^{T}M^{1}(B(Be\cdot Ae))=0,\nonumber\\
&&e^{T}M^{1}(Ae\cdot B^2e)=0,~~e^{T}M^{1}(Ae\cdot (Be)^2)=0,~~e^{T}M^{1}(Be\cdot B A e)=0,\nonumber\\
&&e^{T}M^{1}(B^2e)^2=0,~~e^{T}M^{1}(Be\cdot B^3e)=0,~~e^{T}M^{1}((Be)^2\cdot B^2e)=0,\nonumber\\
&&e^{T}M^{1}(Be\cdot B (Be)^2)=0,~~e^{T}M^{1}(B^4e)=0,~~e^{T}M^{1}(B^2(Be)^2)=0,\nonumber\\
&&e^{T}M^{1}B(B^2e\cdot Be)=0,~~e^{T}M^{1}(B(Be)^3)=0,~~e^{T}M^{1}(BABe)=0,\nonumber\\
&&e^{T}M^{1}(Be\cdot ABe)=0,~~e^{T}M^{1}(A(Be)^2)=0,~~e^{T}M^{1}(AB^2e)=0,\nonumber\\
&&e^{T}M^{1}(Be)^4=0,~~(Be)^{T}M^{1}(Be\cdot B^2e)=0,~~(Be)^{T}M^{1}(B^3e)=0,\nonumber\\
&&(Be)^{T}M^{1}B(Be)^2=0,~~(Be)^{T}M^{1}(Be)^3=0,~~(Be)^{T}M^{1}(BAe)=0,\nonumber\\
&&(Be)^{T}M^{1}(ABe)=0,~~(Be)^{T}M^{1}(Be\cdot Ae)=0,~~(Ae)^{T}M^{1}(Ae)=0,\nonumber\\
&&(Ae)^{T}M^{1}(B^2e)=0,~~(Ae)^{T}M^{1}(Be)^2=0,~~(B^2e)^{T}M^{1}(B^2e)=0,\nonumber\\
&&(B^2e)^{T}M^{1}(Be)^2=0,~~(Be)^{2T}M^{1}(Be)^2=0.
\end{eqnarray*}

Thus, we can obtain a SRK scheme of order $(1.0,2.5-\epsilon)$ as follows
\begin{scheme}\label{scheme2}
\begin{equation}
\begin{tabular}{|cccccccc}
$0$& $0$ & $0$& ~~$0$& $0$ & $0$  \\
$\frac{1}{2}$&  $0$ & $0$&  ~~$\frac{1}{2}$&  $0$ & $0$  \\
$0$&  $1$ &$0$ &~~$0$&  $1$ &$0$  \\
\hline
$\frac{1}{4}$& $\frac{1}{2}$ & $\frac{1}{4}$&~~$\frac{1}{4}$& $\frac{1}{2}$ & $\frac{1}{4}$ \\
\end{tabular}.
\end{equation}
\end{scheme}
\end{example}

\section{Near Preservation of Quadratic Invariants by Implicit SRK Methods}
\label{sec:1}
In the implementation of the implicit SRK method (\ref{SRK-1}), we give a truncation of the random variable $\Delta W_n$ as (see \cite{3})
\begin{eqnarray}\label{d}
Y=D(y_n,Y)=e\otimes y_{n}+h(A\otimes I_d)F(Y)+\overline{\Delta W_n}(B\otimes I_d)G(Y),\nonumber\\
y_{n+1}=y_{n}+h(\alpha^{T}\otimes I_d)F(Y)+\overline{\Delta W_n}(\beta^{T}\otimes I_d)G(Y),~~~~~~~~~~~~~~~~
\end{eqnarray}
for $n=0,1,\cdots,N_T-1$, where, according to \cite{3}, the $\overline{\Delta W_{n}}$ is a truncation of the Wiener increment $\Delta W_{n}$ which satisfies
\begin{equation}\label{truncation}
\overline{\Delta W_n}=\sqrt{h}\xi_h, \quad \xi_h=\left\{\begin{array}{ll}-A_h,&\xi\leq-A_h\\ \xi,&-A_h\leq \xi \leq A_h\\A_h,&\xi\geq A_h\end{array}\right.,
\end{equation}
with $\xi\sim \mathcal{N}(0,1)$, $\Delta W_n=W(t_{n+1})-W(t_n)=\sqrt{h}\xi_h$, and $A_h=\sqrt{2k|\ln h|}$ $(k\geq 1)$. It is proved in \cite{3} that the root-mean-square error of this truncation is $O(h^{\frac{k}{2}})$, thus the mean-square order of the algorithms including such a truncation can be kept the same as that containing the accurate Wiener increment, by appropriate choice of the value $k$.

From the proof of Theorem \ref{thm1} in \cite{15}, it is guaranteed that this truncation will not affect the preservation of quadratic invariants by SRK methods with conditions (\ref{conditionnew}). Therefore, if (\ref{conditionnew}) is satisfied, the implicit SRK method (\ref{d}) preserves quadratic invariants of the underlying SDE. In the following we assume that the conditions in (\ref{conditionnew}) are all satisfied.

Another problem, however, comes from that in implementation, the $Y$ in the first equation of (\ref{d}) can not be solved accurately, but only be approximated by iterations such as fixed-point iteration or Newton's iteration resulting in $Y^{[N]}$ after $N$ iterations. How will the iteration error affect the preservation of the quadratic invariants is what we discuss in the following.
\subsection{Fixed-point iteration}
\label{sec:2}
The fixed-point iteration applied to the first equation of (\ref{d}) takes the form
\begin{eqnarray*}\label{fix}
Y^{[0]}=e\otimes y_n;~~~~~~~~~~~~~~~~~~~~~~~~~~~~~~~~~~~~~~~~~~~~~~~\\
Y^{[N+1]}=D(y_n, Y^{[N]})=e\otimes y_n+hA_1F(Y^{[N]})+\overline{\Delta W_n}B_1G(Y^{[N]}),
\end{eqnarray*}
for $N=0,1,\cdots$, where $A_1=A\otimes I_d$, $B_1=B\otimes I_d$. Now we want to show that, for sufficiently small time step-size $h$, the fixed-point of the mapping $D(y_n,\cdot)$ exists and can really be approximated by performing the iteration (\ref{fix}).

For convenience, we can use the $\|\cdot \|_2$ vector norm as well as its corresponding consistent matrix and tensor norms. In fact, our results hold without dependence on particular choice of norms except for the consistency requirement among them. We use the abbreviation $\parallel\cdot \parallel$ instead of $\|\cdot \|_2$ in the following.

Denote $N(\Omega,\epsilon)=\{y\in \mathbb{R}^d: \underset{y_n\in{\Omega}}{\min}  \parallel y-y_n \parallel\leq \epsilon\}$, and $N^s(\Omega,\epsilon)=\underbrace{N(\Omega,\epsilon)\times\cdots \times N(\Omega,\epsilon)}_{s}$. We prove a lemma similar to Proposition 1 in \cite{paper}. Note that all the results we obtained  hold in the sense of `almost surely'.

\begin{lemma}\label{fixed thm} Let $\Omega\subset\mathcal{R}^d$ be a bounded, convex and open set, and $f$, $g$ be globally Lipschitz continuous on $\Omega$ with Lipschitz constants $L$ and $M$, respectively. Then, for any $\epsilon>0$, there exists $h_0>0$ dependent on $\Omega$ and $\epsilon$, such that for any $h\leq h_0$, and $y_n\in \Omega$,
\begin{enumerate}
\item $D(y_n,\cdot)$ maps $N^s(\Omega,\epsilon)$ into itself;
\item There exists a unique solution $Y^{*}$ to the first equation of (\ref{d}), and it can be approximated via the iteration (\ref{fix});
\item $\| Y^{[N]}-Y^{*}\|\leq \delta^N \| Y^0-Y^* \|$ with $0<\delta<1$, where $\delta=C_1\sqrt{2kh|\ln h|}$ and $C_1=\max\{\|A_1\|L, \|B_1\|M\}$.
\end{enumerate}
\end{lemma}
$Proof.$
$\forall Y\in N^s(\Omega,\epsilon),$ $\forall h\leq e^{-\frac{1}{2k}}$,
\begin{eqnarray*}\label{proof1}
\|D(y_n,Y)-e\otimes y_n\|&\leq &h \|A_1\|C_0+\|\overline{\Delta W_n}\|\|B_1\|\tilde{C}_0\\
&\leq&\sqrt{2kh|\ln h|}( \|A_1\|C_0+\|B_1\|\tilde{C}_0),
\end{eqnarray*}
where $C_0=\underset{Y\in N^s(\Omega,\epsilon)}{\max}\|F(Y)\|$, $\tilde{C}_0=\underset{Y\in N^s(\Omega,\epsilon)}{\max}\|G(Y)\|$. Since $\sqrt{h|\ln h|}\rightarrow 0$ as $h\rightarrow 0$, $\forall \epsilon >0$, there exists $0<h_1\leq e^{-\frac{1}{2k}}$, such that for all $h\leq h_1$, and $\forall y_n\in \Omega$, $\|D(y_n,Y)-e\otimes y_n\|< \epsilon.$ Thus (a) is verified. Next we prove that $D(y_n,\cdot)$ is a contraction mapping for sufficiently small $h$. $\forall Y^1,\,\,Y^2\in N^s(\Omega,\epsilon)$,
\begin{eqnarray*}\label{proof2}
\|D(y_n,Y^1)-D(y_n,Y^2)\|&=\|hA_1(F(Y^1)-F(Y^2))+\overline{\Delta W_n}B_1(G(Y^1)-G(Y^2))\|\\
&\leq \sqrt{2kh|\ln h|}(\|A_1\|L+\|B_1\|M)\|Y^1-Y^2\|.
\end{eqnarray*}
There exists $0<h_2\leq e^{-\frac{1}{2k}}$, such that $\forall h\leq h_2$, $\forall y_n\in \Omega$,
\begin{equation}\label{1}
\sqrt{2kh|\ln h|}(\|A_1\|L+\|B_1\|M)<1,
\end{equation}
which implies that $D(y_n,\cdot)$ is a contraction mapping. Let $h_0=\min\{h_1,h_2\}$. Then, by the contraction mapping principle, $\forall h<h_0$, there exists a unique solution $Y^*$ to the first equation of (\ref{d}), which can be approximated via the iteration (\ref{fix}), and
\begin{eqnarray*}\label{proof3}
\|Y^{[N]}-Y^*\|&=&\|hA_1(F(Y^{[N-1]})-F(Y^*))+\overline{\Delta W_n}B_1(G(Y^{[N-1]})-G(Y^*))\|\\
&\leq& C_1\sqrt{2kh|\ln h|}\|Y^{[N-1]}-Y^*\|\leq\cdots\\
&\leq& \delta^N\|Y^{[0]}-Y^*\|.
\end{eqnarray*}
It follows immediately from (\ref{1}) that $0<\delta<1$. This completes the proof.

Denote
\begin{equation}\label{yN}
y_{n+1}^{[N]}=y_{n}+h(\alpha^T\otimes I_d)F(Y^{[N]})+\overline{\Delta W_n}(\beta^T\otimes I_d)G(Y^{[N]}).
\end{equation}
Now we estimate the error in the preservation of the quadratic invariants $y^T C y$ caused by the iterations (\ref{fix}) on $Y$. We have the following result.
\begin{theorem}\label{fixed thm2} Let $\Omega$, $f$ and $g$ satisfy the same assumptions as in Lemma \ref{fixed thm}, and the conditions (\ref{conditionnew}) hold. Let $h_0$ be the one given in Lemma \ref{fixed thm}. Then $\forall h\leq h_0$, $\forall y_n\in \Omega$,
\begin{equation}\label{thm eq}
|(y_{n+1}^{[N]})^T C y_{n+1}^{[N]}-y_n^TCy_n|\leq\|C\|[\frac{C_2^2D_1^2}{C_1^4}\delta ^{2N+4}+\frac{2C_2D_0D_1}{C_1^2}\delta^{N+2}],
\end{equation}
where $C_2=\|\alpha\|L+\|\ \beta\|M$, $D_0=\|y_0\|+R_0$, with $R_0$ being the diameter of $\Omega$, and $D_1=\|A_1\|C_0+\|B_1\|\tilde{C}_0$.
\end{theorem}

$Proof.$
Under the conditions (\ref{conditionnew}), the implicit SRK method (\ref{d}) preserves the quadratic invariant, i.e.,
\begin{equation}\label{preservation}
y_{n+1}^TCy_{n+1}-y_n^TCy_n=0.
\end{equation}
Therefore,
\begin{eqnarray*}\label{difference}
&\quad\,\,(y_{n+1}^{[N]})^TCy_{n+1}^{[N]}-y_n^TCy_n=(y_{n+1}^{[N]})^TCy_{n+1}^{[N]}-y_{n+1}^TCy_{n+1}\\
&=(y_{n+1}^{[N]}-y_{n+1})^TC(y_{n+1}^{[N]}-y_{n+1})+2(y_{n+1}^{[N]}-y_{n+1})^TCy_{n+1}.
\end{eqnarray*}
Thus,
\begin{equation}\label{proof22}
|(y_{n+1}^{[N]})^TCy_{n+1}^{[N]}-y_n^TCy_n|\leq \|C\|[\|y_{n+1}^{[N]}-y_{n+1}\|^2+2\|y_{n+1}^{[N]}-y_{n+1}\|\|y_{n+1}\|].
\end{equation}
Meanwhile, for $h\leq h_0$, we have
\begin{equation}\label{proof23}
\|y_{n+1}^{[N]}-y_{n+1}\|\leq(h\|\alpha\|L+\sqrt{2kh|\ln h|}\|\beta\|M)\|Y^{[N]}-Y^*\|\\
\leq \frac{\delta}{C_1}C_2\delta^N\|Y^{[0]}-Y^*\|,
\end{equation}
and
\begin{eqnarray}\label{proof24}
\|Y^{[0]}-Y^*\|&\leq h\|A_1\|C_0+\sqrt{2kh|\ln h|}\|B_1 \|\tilde{C}_0\leq\frac{\delta}{C_1}D_1,\nonumber\\
&\|y_{n+1}\|\leq \|y_0\|+R_0.
\end{eqnarray}
Substitute (\ref{proof23}) and (\ref{proof24}) into (\ref{proof22}), we derive the result (\ref{thm eq}). This completes the proof.

\begin{remark}\label{rm4}
Since $y_n^TCy_n=y_0^TCy_0$, the inequality (\ref{thm eq}) implies that, there exist constants $K_1$ and $K_2$ depending on $N_T$, such that
\begin{equation}\label{thm eq-1}
|I(y_{N_T}^{[N]})-I(y_0)|\leq K_1\delta ^{2N+4}+K_2\delta^{N+2}.
\end{equation}
\end{remark}
\subsection{Newton's iteration}
\label{sec:2}
Newton's iteration applied to the first equation of (\ref{d}) reads
\begin{eqnarray}\label{newt}
&Y^{[0]}=e\otimes y_n,\nonumber\\
&Y^{[N+1]}=\hat{D}(y_n,Y^{[N]})\nonumber\\
&=Y^{[N]}-[I_{sd}-hA_1F'(Y^{[N]})-\overline{\Delta W_n}B_1G'(Y^{[N]})]^{-1}\nonumber\\&\cdot[Y^{[N]}-e\otimes y_n-hA_1F(Y^{[N]})-\overline{\Delta W_n}B_1G(Y^{[N]})].
\end{eqnarray}
Choose $0<\tau_1<e^{-\frac{1}{2k}}$ such that $\forall h\leq \tau_1$, $\forall Y\in N^s(\Omega,\epsilon)$, the matrix $I_{sd}-hA_1F'(Y)-\overline{\Delta W_n}B_1G'(Y)$ is invertible. We establish the following lemma which is a stochastic extension of the corresponding deterministic result in \cite{paper}.
\begin{lemma}\label{lemma newt}
Let $\Omega$ be assumed as in Lemma \ref{fixed thm}, and $f$, $g$ be three times continuously differentiable on $\Omega$. Then for any $\epsilon>0$, there exists $\tau_0>0$ dependent on $\Omega$ and $\epsilon$, such that $\forall h\leq \tau_0$, $\forall y_n\in \Omega$,

\begin{enumerate}
\item  $\hat{D}(y_n,\cdot)$ maps $N^s(\Omega,\epsilon)$ into itself;
\item There exists a unique solution $Y^{*}$ to the first equation of (\ref{d}), and it can be approximated by the iteration (\ref{newt});
\item $\| Y^{[N]}-Y^{*}\|\leq \gamma^{2^N-1} \| Y^{[0]}-Y^* \|^{2^N}$ with $\gamma>0$, and $\gamma\ \cdot| Y^{[0]}-Y^* \|<1$.
\end{enumerate}
\end{lemma}
$Proof.$ The proof is similar to that in \cite{paper}, but extended to stochastic context.

Denote $R=hA_1F'(Y)+\overline{\Delta W_n}B_1G'(Y)$. Since
\begin{eqnarray*}
(I_{sd}-R)^{-1}=I_{sd}+(I_{sd}-R)^{-1}R,
\end{eqnarray*}
we have
\begin{align}\label{proof30}
\hat{D}&(y_n,Y)\nonumber\\
&=Y-[I_{sd}+(I_{sd}-R)^{-1}R][Y-e\otimes y_n-(hA_1F(Y)+\overline{\Delta W_n}B_1G(Y))]\nonumber\\
&=e\otimes y_n+(I_{sd}-R)^{-1}[R(e\otimes y_n-Y)+hA_1F(Y)+\overline{\Delta W_n}B_1G(Y)].
\end{align}
Denote
\begin{eqnarray*}
C_3=\underset{Y\in N^s(\Omega,\epsilon)}{\max}\|F'(Y)\|,\quad \tilde{C}_3=\underset{Y\in N^s(\Omega,\epsilon)}{\max}\|G'(Y)\|,\\
D_2=\|A_1\|C_3+\|B_1\|\tilde{C}_3,\quad D_3=\underset{Y\in N^s(\Omega,\epsilon),y_n\in \Omega}{\max}\|e\otimes y_n-Y\|,
\end{eqnarray*}
then for $h\leq \tau_1$,
\begin{eqnarray}\label{proof31}
\|R\|\leq \sqrt{2kh|\ln h|}D_2,~~~~~~~~~~~~~~~~~\nonumber\\
 \|(I_{sd}-R)^{-1}\|\leq \frac{1}{1-\|R\|}\leq\frac{1}{1-\sqrt{2kh|\ln h|}D_2}\nonumber\\
\|hA_1F(Y)+\overline{\Delta W_n}B_1G(Y)\|\leq \sqrt{2kh|\ln h|}D_1.
\end{eqnarray}
Substitute (\ref{proof31}) into (\ref{proof30}), we have for $h\leq \tau_1$,
\begin{equation}\label{proof32}
\|\hat{D}(y_n,Y)-e\otimes y_n\|\leq \frac{\sqrt{2kh|\ln h|}}{1-D_2\sqrt{2kh|\ln h|}}(D_2D_3+D_1).
\end{equation}
Since $\frac{\sqrt{2kh|\ln h|}}{1-D_2\sqrt{2kh|\ln h|}}\rightarrow 0$ as $h\rightarrow 0$, $\forall \epsilon>0$, $\forall y_n\in \Omega$, there exists $\tau_2>0$ depending on $\epsilon$ and $\Omega$, such that $\forall h\leq \tau_2$,
\begin{equation}\label{proof33}
\frac{\sqrt{2kh|\ln h|}}{1-D_2\sqrt{2kh|\ln h|}}(D_2D_3+D_1)<\epsilon.
\end{equation}
Choose $\tau_3=\min\{\tau_1,\tau_2\}$, then for $h\leq \tau_3$,
$\forall y_{n}\in \Omega$, $\hat{D}(y_n,\cdot)$ maps
$N^s(\Omega,\epsilon)$ into itself. Thus (a) is proved.

To prove $\hat{D}$ is a contraction mapping, we observe
$\frac{\partial \hat{D}}{\partial Y}(y_n,Y)$. Since
$\frac{\partial}{\partial Y}(I_{sd}-R)^{-1}=(I_{sd}-R)^{-2}R'(Y)$, it holds
for any vector $\eta\in \mathcal{R}^{sd}$,
\begin{eqnarray*}\label{tensor}
\frac{\partial \hat{D}}{\partial
Y}(y_n,Y)\eta&=-(I_{sd}-R)^{-1}[hA_1F''(Y)\eta+\overline{\Delta
W_n}B_1G''(Y)\eta](I_{sd}-R)^{-1}\\&[Y-e\otimes
y_n-hA_1F(Y)-\overline{\Delta W_n}B_1G(Y)].
\end{eqnarray*}
Denote
\begin{eqnarray*}\label{c4}
C_4=\underset{Y\in N^s(\Omega,\epsilon)}{\max}\|F''(Y)\|,\quad
\tilde{C}_4=\underset{Y\in
N^s(\Omega,\epsilon)}{\max}\|G''(Y)\|\\
D_4=\|A_1\|C_4+\|B_1\|\tilde{C}_4.~~~~~~~~~~~~~~~~~~~
\end{eqnarray*}
$\forall h\leq \tau_1$, according to (\ref{proof31}), we have
\begin{eqnarray}\label{partial}
\|\frac{\partial \hat{D}}{\partial Y}\|\leq
\frac{D_3D_4\sqrt{2kh|\ln h|}+D_1D_4(2kh|\ln
h|)}{(1-D_2\sqrt{2kh|\ln h|})^2}.
\end{eqnarray}
It is obvious that the right hand side of (\ref{partial}) tends to
zero as $h\rightarrow 0$. Therefore, there exists
$0<\tau_4\leq\tau_1$, such that $\forall h\leq \tau_4$, $\forall
y_n\in \Omega$
\begin{equation}\label{t4}
\|\frac{\partial \hat{D}}{\partial Y}\|<1,
\end{equation}
which implies that $\hat{D}$ is a contraction mapping. Let
$\tau_5=\min\{\tau_3,\tau_4\}$. Then, according to the contraction
mapping principle, for all $h\leq\tau_5$, $\forall y_n\in \Omega$,
there exists a unique solution $Y^*$ to the first equation of
(\ref{d}), which can be approximated via the iteration (\ref{newt}).

Equation (\ref{tensor}) implies $\frac{\partial\hat{D}}{\partial
Y}(y_n,Y^*)=\mathbf{0}$. Therefore the iteration (\ref{newt})
converges in the second order, i.e.,
\begin{eqnarray*}\label{newt diff}
\|Y^{[N]}-Y^*\|&\leq\gamma \|Y^{[N-1]}-Y^*\|^2\leq\cdots\\
&\leq \gamma^{2^N-1}\|Y^{[0]}-Y^*\|^{2^N},
\end{eqnarray*}
where $\gamma=\underset{Y\in N^s(\Omega,\epsilon),y_n\in\Omega,
h\leq \tau_5}{\max}\|\frac{\partial^2 \hat{D}}{\partial
Y^2}(y_n,Y)\|$.
It is not difficult to check that $\gamma\rightarrow
0$ as $h\rightarrow 0$. On the other hand, $\|Y^*-Y^{[0]}\|\leq
\sqrt{2kh|\ln h|}D_1$. Consequently, there exists $\tau_6>0$ such
that for $h\leq \tau_6$, $\gamma\cdot\|Y^*-Y^{[0]}\|<1 $. Choose
$\tau_0=\min\{\tau_5,\tau_6\}$, the lemma is fully verified. This
completes the proof.

Let $Y^{[N]}$ result from $N$ times iteration via (\ref{newt}), and
$y_{n+1}^{[N]}$ be defined as in (\ref{yN}). Denote
\begin{equation}\label{delta1}
\hat{\delta}=\gamma D_1\sqrt{2kh|\ln h|}.
\end{equation}
Then there exists $\tau_7>0$, such that for $h\leq \tau_7$,
$0<\hat{\delta}<1$. Let $\hat{\tau}_0=\min\{\tau_0,\tau_7\}$. We
have the following theorem estimating the error in preservation of
the quadratic invariants caused by iterations via (\ref{newt}).

\begin{theorem}\label{ntthm}
Let $\Omega$, $f$ and $g$ satisfy the same assumptions as in Lemma \ref{lemma newt}, and the conditions (\ref{conditionnew}) hold. Then $\forall h\leq
\hat{\tau}_0$, $\forall y_n\in \Omega$,
\begin{equation}\label{thm1 eq}
|(y_{n+1}^{[N]})^T C
y_{n+1}^{[N]}-y_n^TCy_n|\leq\|C\|[\frac{C_2^2}{D_1^2\gamma^4}\hat{\delta}
^{2^{N+1}+2}+\frac{2C_2D_0}{D_1\gamma^2}\hat{\delta}^{2^{N}+1}],
\end{equation}
with $0<\hat{\delta}<1$.
\end{theorem}

$Proof.$ For $h\leq \hat{\tau}_0$, using the results
(\ref{proof24}) and (\ref{newt diff}), we have
\begin{eqnarray*}\label{proof4}
\|y_{n+1}^{[N]}-y_{n+1}\|&\leq& [h\|\alpha\|L+\sqrt{2kh|\ln
h|}\|\beta\|M]\|Y^{[N]}-Y^*\|\\&\leq& \sqrt{2kh|\ln h|}C_2
\gamma^{2^N-1}\|Y^{[0]}-Y^*\|^{2^N}\\&\leq&
\frac{C_2}{D_1\gamma^2}\hat{\delta}^{2^N+1}.
\end{eqnarray*}
Similar to (\ref{difference}),
\begin{eqnarray*}\label{difference newt}
|(y_{n+1}^{[N]})^TCy_{n+1}^{[N]}-y_n^TCy_n|&\leq
\|C\|[\|y_{n+1}^{[N]}-y_{n+1}\|^2+2\|y_{n+1}^{[N]}-y_{n+1}\|\|y_{n+1}\|]\\&\leq\|C\|[\frac{C_2^2}{D_1^2\gamma^4}\hat{\delta}
^{2^{N+1}+2}+\frac{2C_2D_0}{D_1\gamma^2}\hat{\delta}^{2^{N}+1}].
\end{eqnarray*}
This completes the proof.
\begin{remark}\hfill
\begin{itemize}
\item Similar to Remark \ref{rm4}, (\ref{thm1 eq}) implies that, there exist constant $\hat{K}_1$ and $\hat{K}_2$ depending on $N_T$ such that
\begin{equation}\label{thm1 eq-1}
|I(y_{N_T}^{[N]})-I(y_0)|\leq\hat{K}_1\hat{\delta}
^{2^{N+1}+2}+\hat{K}_2\hat{\delta}^{2^{N}+1}.
\end{equation}
\item
As discussed in Remark \ref{rm1}, the implicit SRK method (\ref{d}) with
conditions (\ref{conditionnew}) applied to (\ref{SHS}) can preserve the
symplectic structure of (\ref{SHS}) accurately, though with truncation of $\Delta W_n$ (see \cite{3}). In implementation of the implicit
SRK methods (\ref{d}), the error in the
preservation of the symplectic structure
$$\|(\psi^{[N]}_{n+1})^TJ\psi^{[N]}_{n+1}-\psi_{n}^TJ\psi_n\|$$ arising
from fixed-point or Newton's iterations can be estimated according
to the results in Theorem \ref{fixed thm2} and \ref{ntthm}, respectively.
\item Compared to the results on approximate preservation of symplectic structure by deterministic symplectic
Runge-Kutta methods due to iterations in implementation (see \cite{paper}), we find that for
fixed-point iteration, the leading error term in the symplectic structure by deterministic symplectic
Runge-Kutta methods is $O(h^{N+2})$, while that by
stochastic symplectic
Runge-Kutta methods is $O((h|\ln h|)^\frac{N+2}{2})$. For
Newton's iteration, the leading error term in the symplectic structure by deterministic and stochastic symplectic
Runge-Kutta methods are $O(h^{2^N-1})$ and
$O((h|\ln h|)^{2^{N-1}+\frac{1}{2}})$, respectively.
\end{itemize}
\end{remark}
\section{Numerical Experiments}
\label{sec;preliminary}
In this section, we apply the explicit SRK schemes, Scheme \ref{scheme1} and \ref{scheme2}, as well as the the stochastic midpoint rule
to test the behavior of the schemes in preserving quadratic invariants,
implemented with fixed point iteration to the Kubo oscillator, and with Newton's iteration to a non-linear stochastic Hamiltonian system.
\subsection{Explicit SRK methods}
\label{sec:2}
The Kubo oscillator
\begin{eqnarray*}\label{kubo}
dp&=-aqdt-\sigma q\circ dW(t),\quad p(0)=p_0,\\dx&=apdt+\sigma p\circ dW(t),\qquad x(0)=x_0,
\end{eqnarray*}
where $a$ and $\sigma$ are constants, is a stochastic system with the quadratic invariant (see \cite{3})
\begin{equation}\label{circle}
H(p(t),x(t))=p(t)^2+x(t)^2=p_0^2+x_0^2,\quad \forall t\geq 0.
\end{equation}
That is, the phase trajectory is a circle with center at the origin and radius $p_0^2+x_0^2$.

\begin{figure}[htbp]
\begin{center}
\includegraphics[height=5cm,width=5cm]{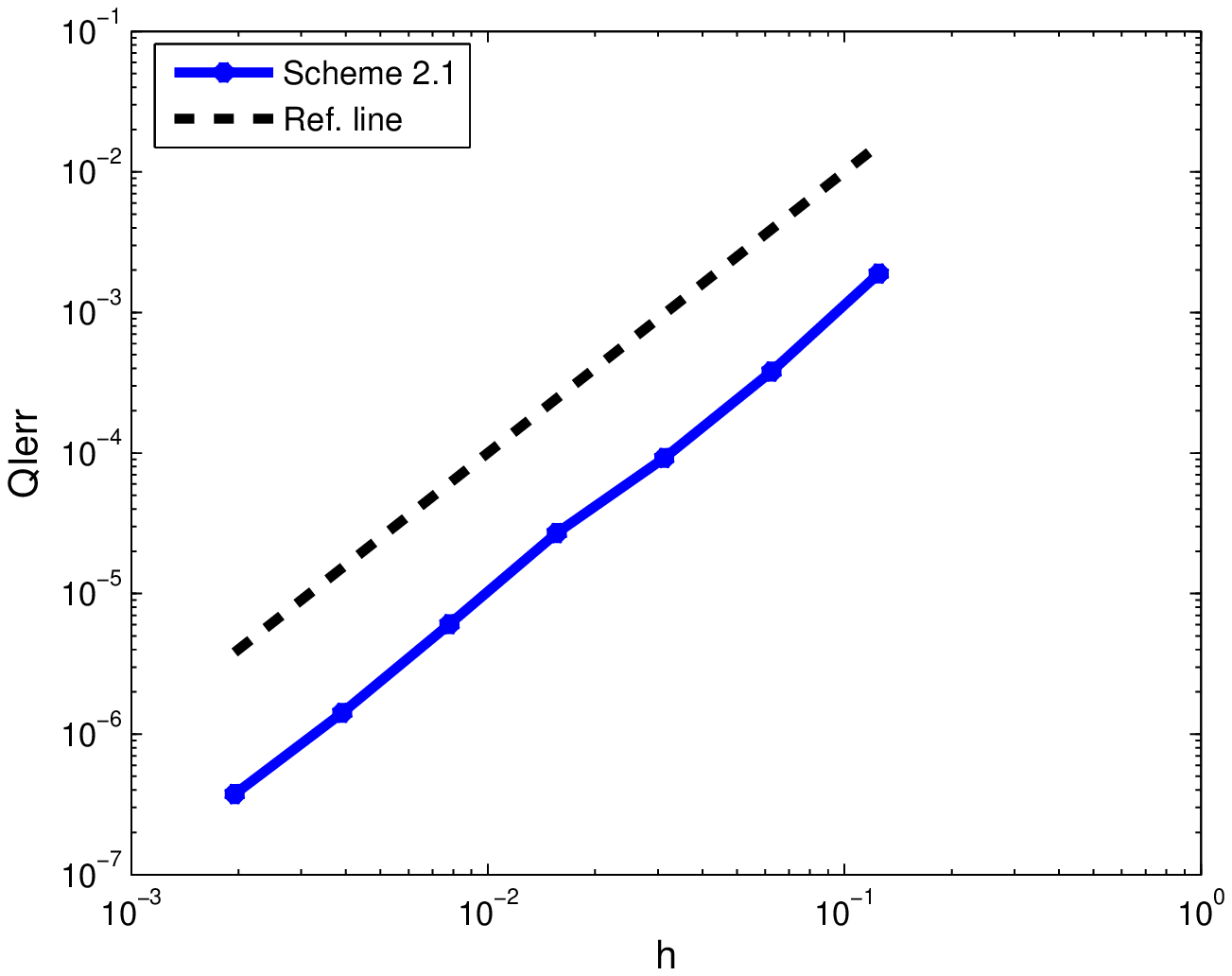}
\includegraphics[height=5cm,width=5cm]{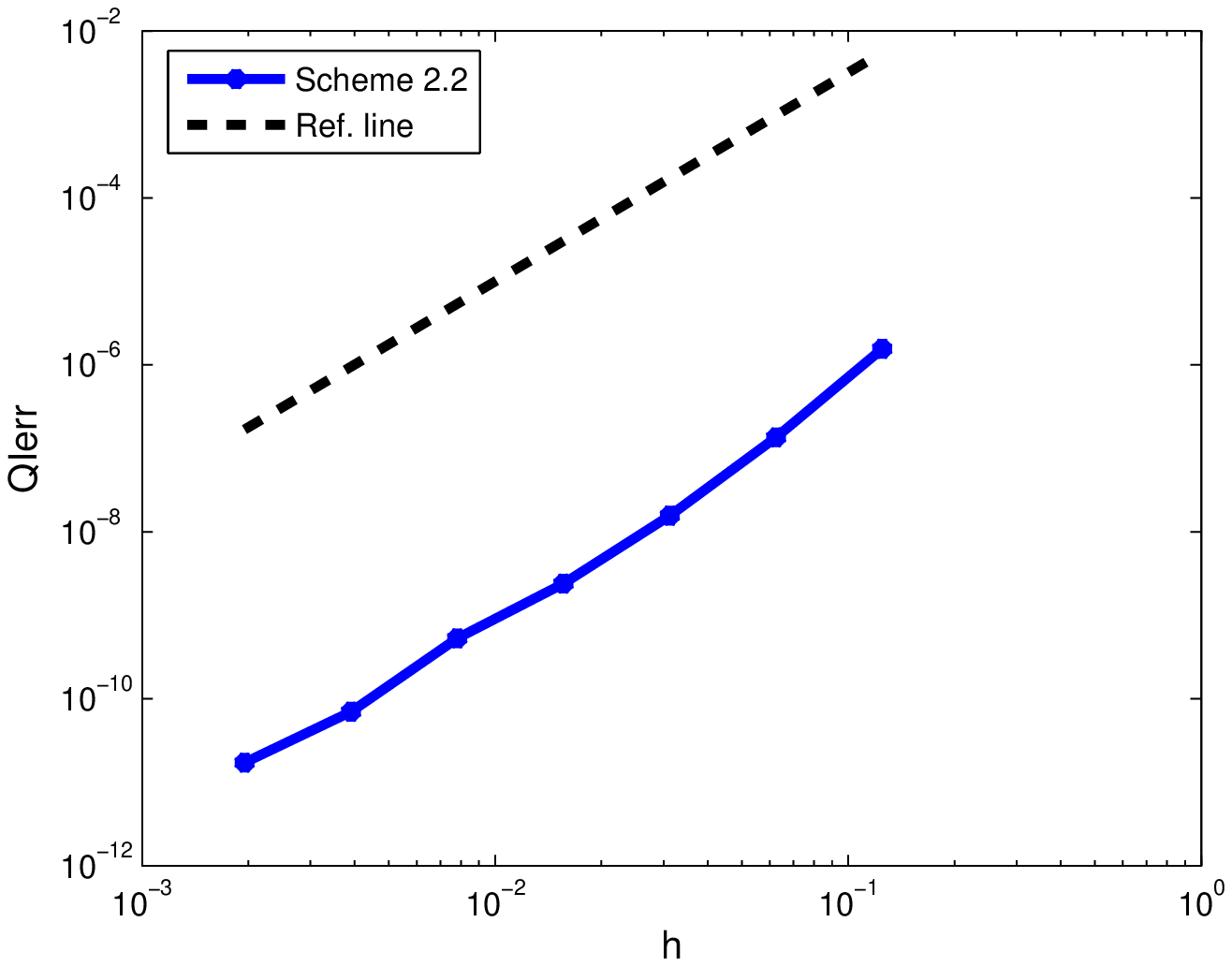}
\caption{The error of preservation of $H(p,x)$ produced by the Scheme 2.1  versus step-size (left), and by Scheme 2.2  versus step-size (right), respectively. The reference line is of slope $2.0$ (left) and $2.5$ (right), respectively.}
\label{figure01}
\end{center}
\end{figure}

We apply Schemes \ref{scheme1} and \ref{scheme2} to this system to observe the order of preservation of quadratic invariants of these two schemes. Choose $a=1$, $\sigma=1$, $T=1$ and $(p_0,x_0)=(0,1)$. Figure \ref{figure01} shows that the order of preservation of quadratic invariants of these two schemes are in good accordance with the theoretically predicted order  $2.0-\epsilon$ and $2.5-\epsilon$, respectively.

Let $a=1$. Now we compare the Schemes \ref{scheme1} and \ref{scheme2} in preserving quadratic invariants $H(p,x)$ with the Milstein scheme (MI-scheme) \begin{equation}\label{milsteinscheme}
\begin{split}
x_{n+1}&=x_{n}+p_{n}h+\sigma p_n \Delta W_n-\frac{1}{2}\sigma^2 x_{n}\Delta W_n^2,\\
p_{n+1}&=p_{n}-x_{n}h-\sigma x_n\Delta W_n-\frac{1}{2}\sigma^2 p_{n}\Delta W_n^2.
\end{split}
\end{equation}
and the midpoint scheme (Midscheme)
\begin{equation}\label{midscheme}
\begin{split}
x_{n+1}&=x_{n}+h\frac{p_{n}+p_{n+1}}{2}+\sigma \Delta W_n\frac{p_{n}+p_{n+1}}{2},\\
p_{n+1}&=p_{n}-h\frac{x_{n}+x_{n+1}}{2}-\sigma \Delta W_n\frac{x_{n}+x_{n+1}}{2},
\end{split}
\end{equation}
which are both of convergence order $1.0$ as the Schemes \ref{scheme1} and \ref{scheme2}. Note that the implicit midpoint scheme reduces in this example to an explicit one due to linearity of the system (\ref{kubo}). The Butcher tabular for the stochastic midpoint scheme is
\begin{tabular}{c|cc}&$\frac{1}{2}$&$\frac{1}{2}$\\ \hline
&1&1
\end{tabular}, which implies that it satisfies the condition (\ref{conditionnew}) of preserving quadratic invariants.

\begin{figure}[htbp]
\begin{center}
\includegraphics[height=5cm,width=5cm]{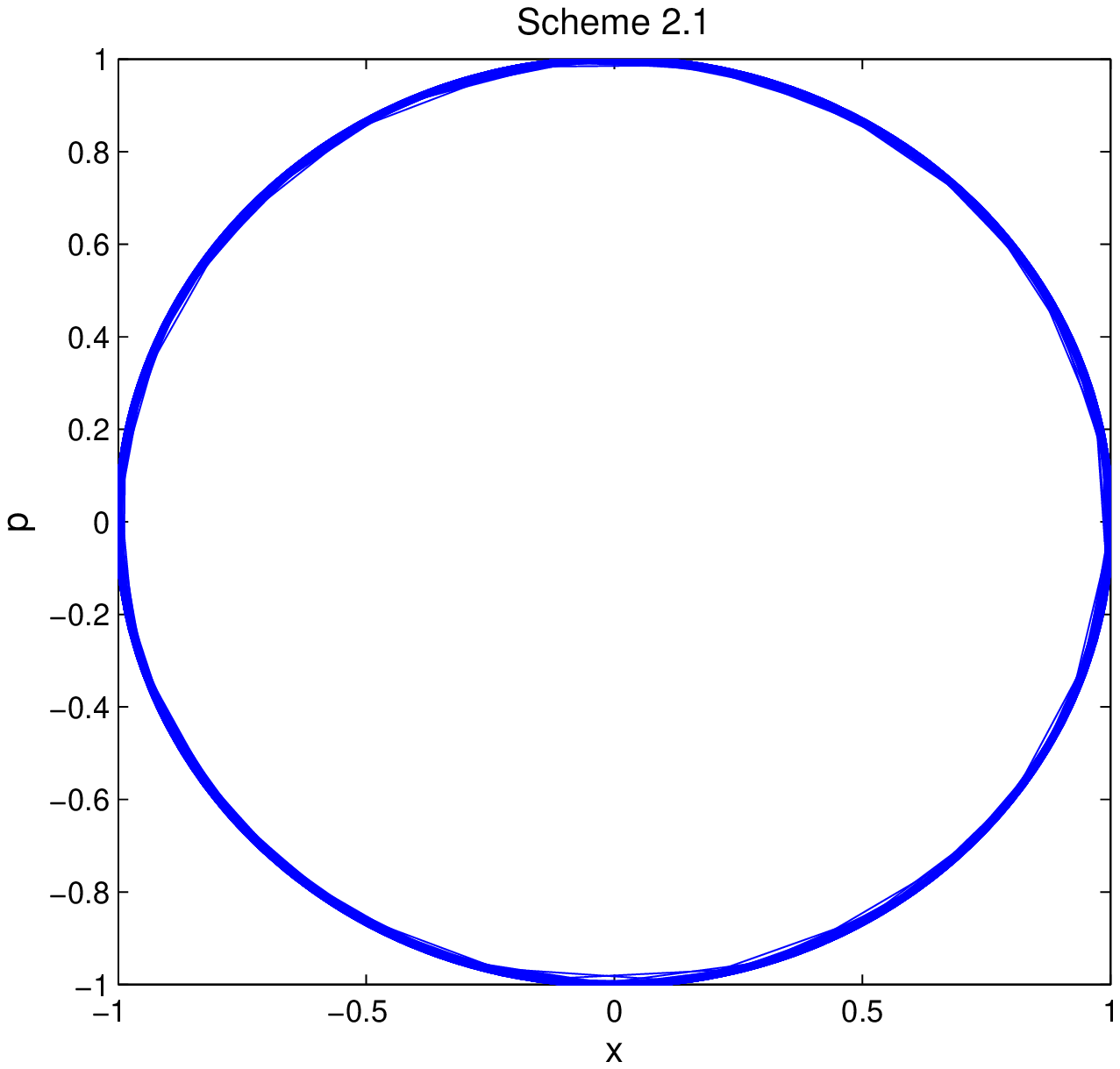}
\includegraphics[height=5cm,width=5cm]{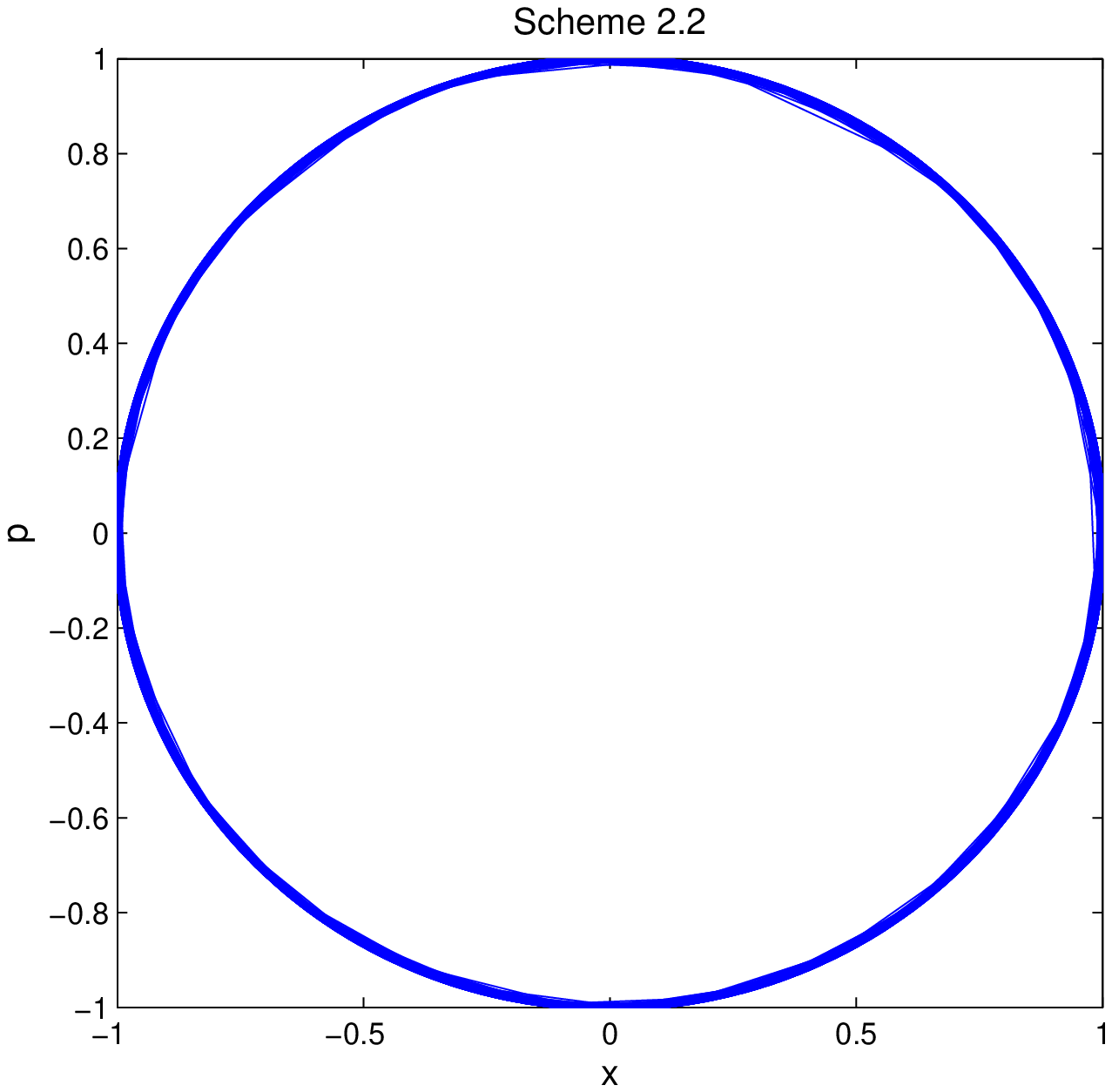}\\
\includegraphics[height=5cm,width=5cm]{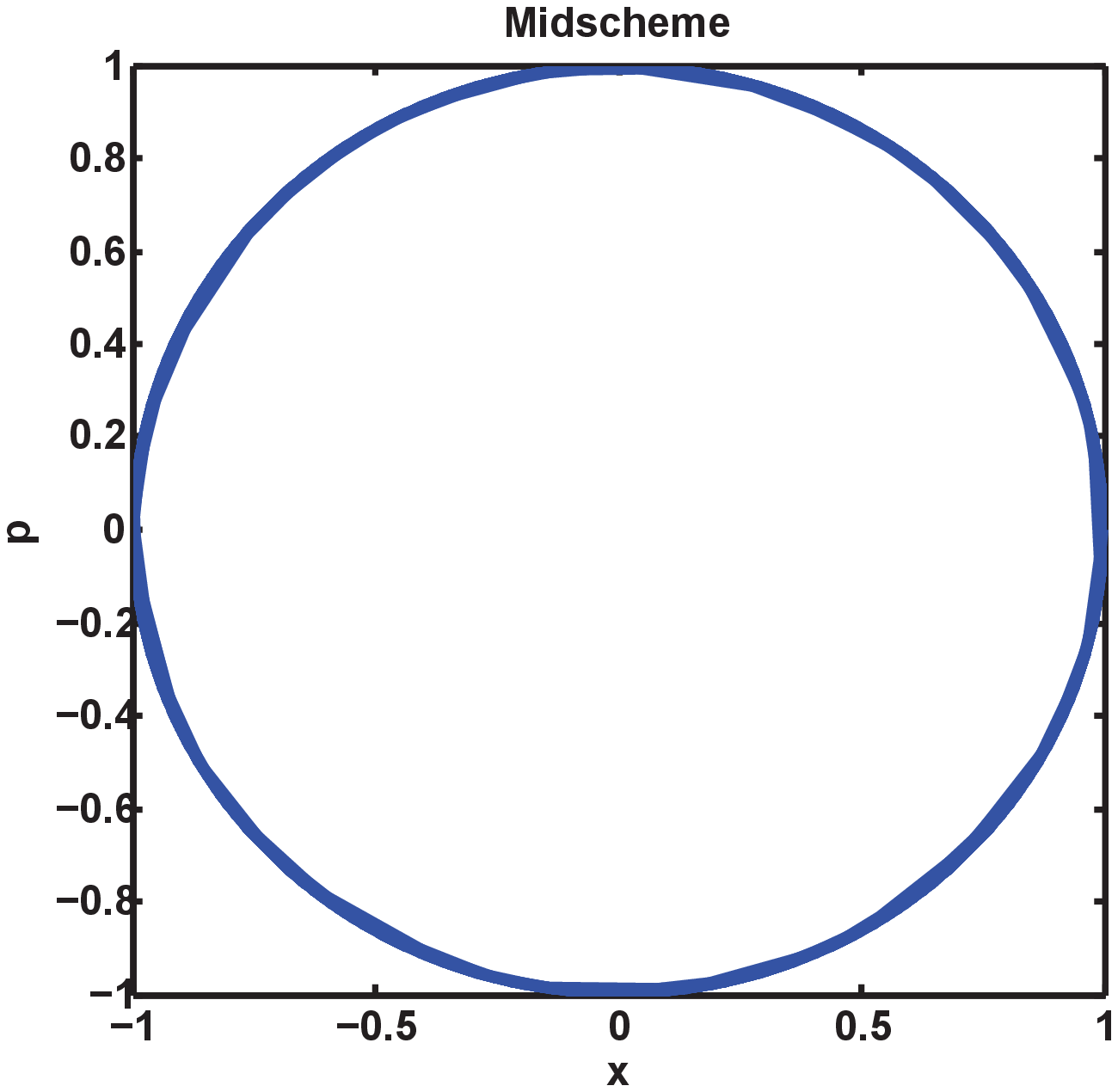}
\includegraphics[height=5cm,width=5cm]{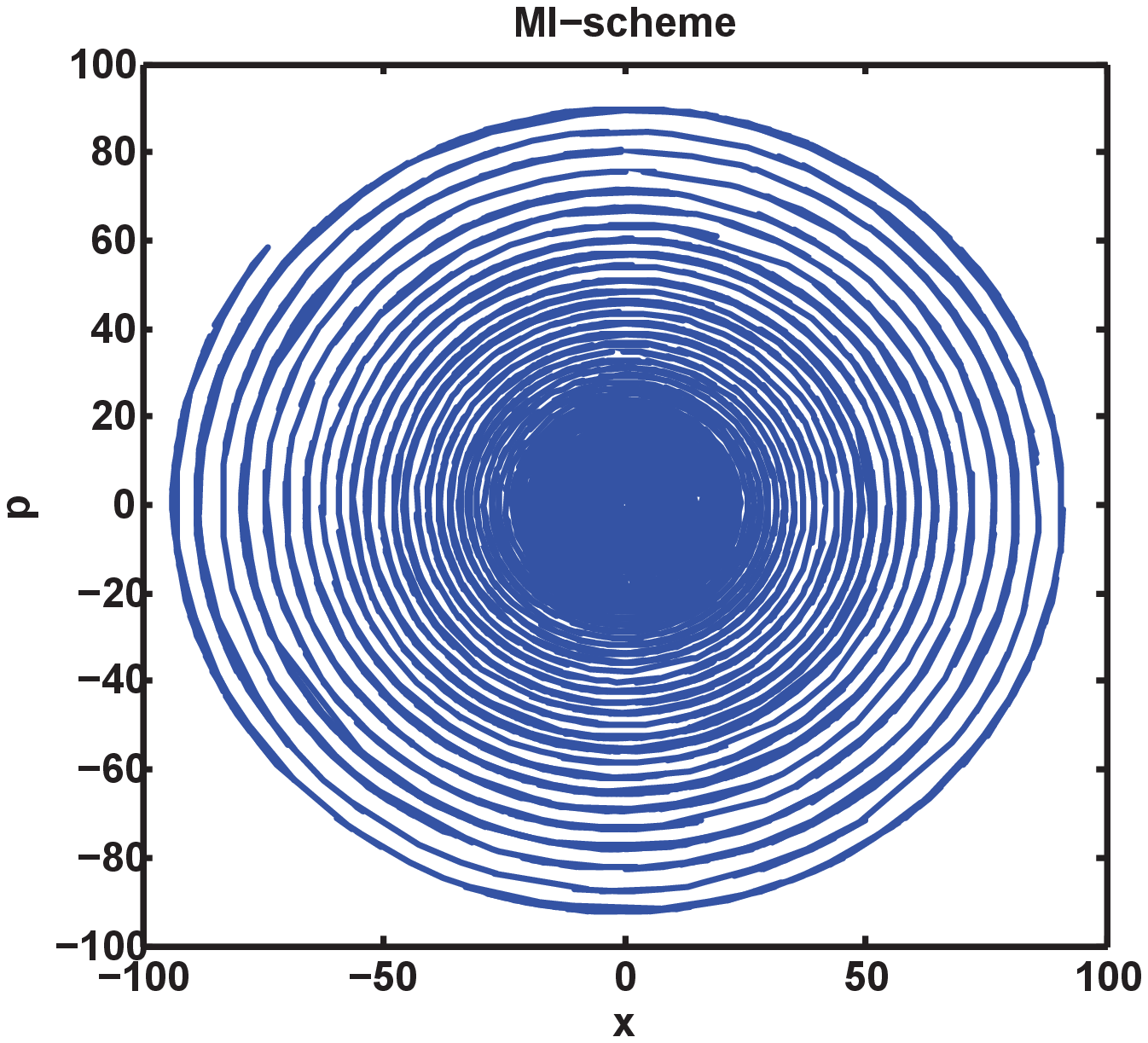}
\caption{Numerical phase trajectory produced by Scheme 2.1 (top left), Scheme 2.2 (top right), midpoint scheme (bottom left) and  Milstein scheme (bottom right) for $\sigma=1$.}
\label{figure03}
\end{center}
\end{figure}

\begin{figure}[htbp]
\begin{center}
\includegraphics[height=5cm,width=5cm]{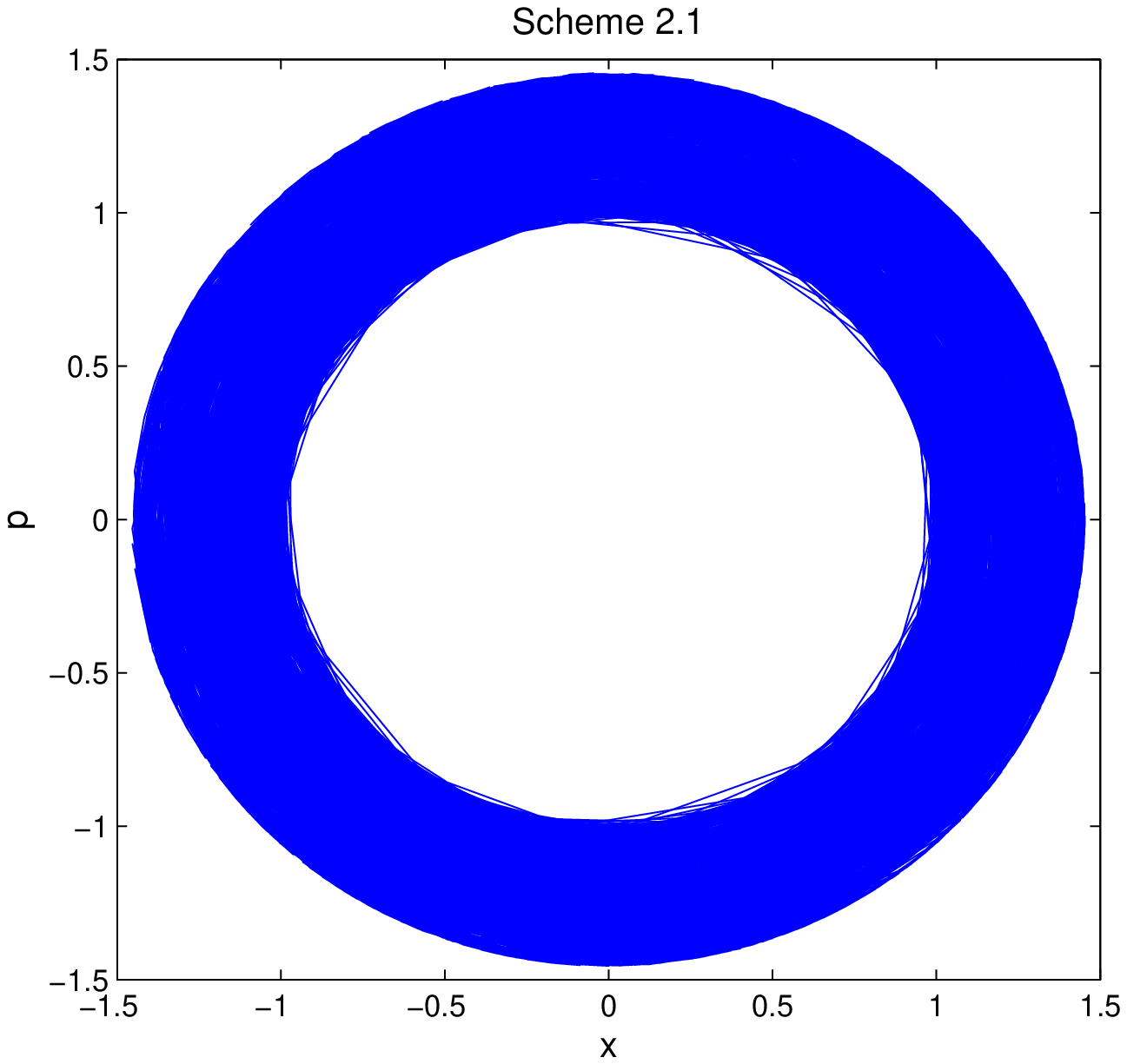}
\includegraphics[height=5cm,width=5cm]{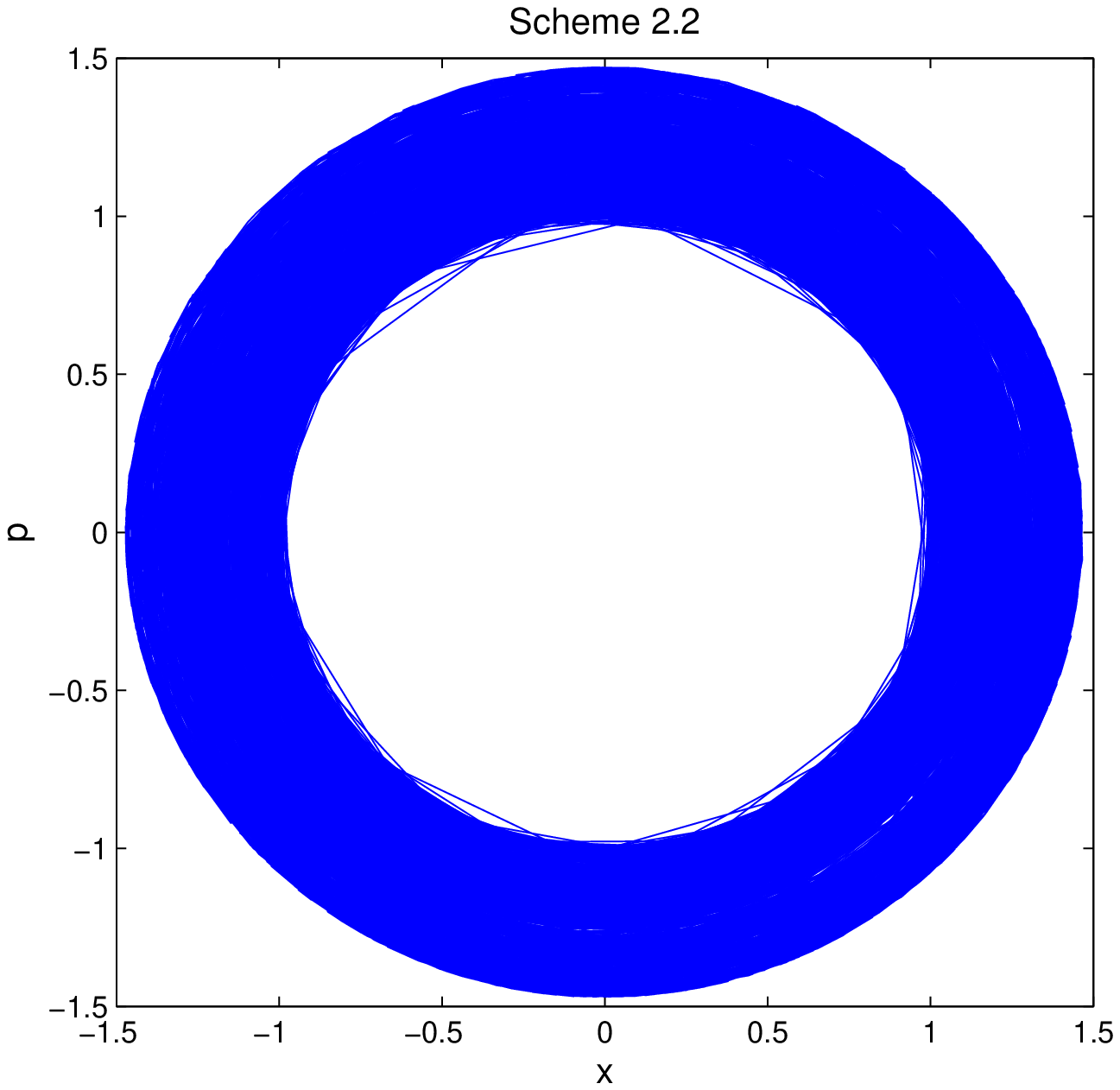}\\
\includegraphics[height=5cm,width=5cm]{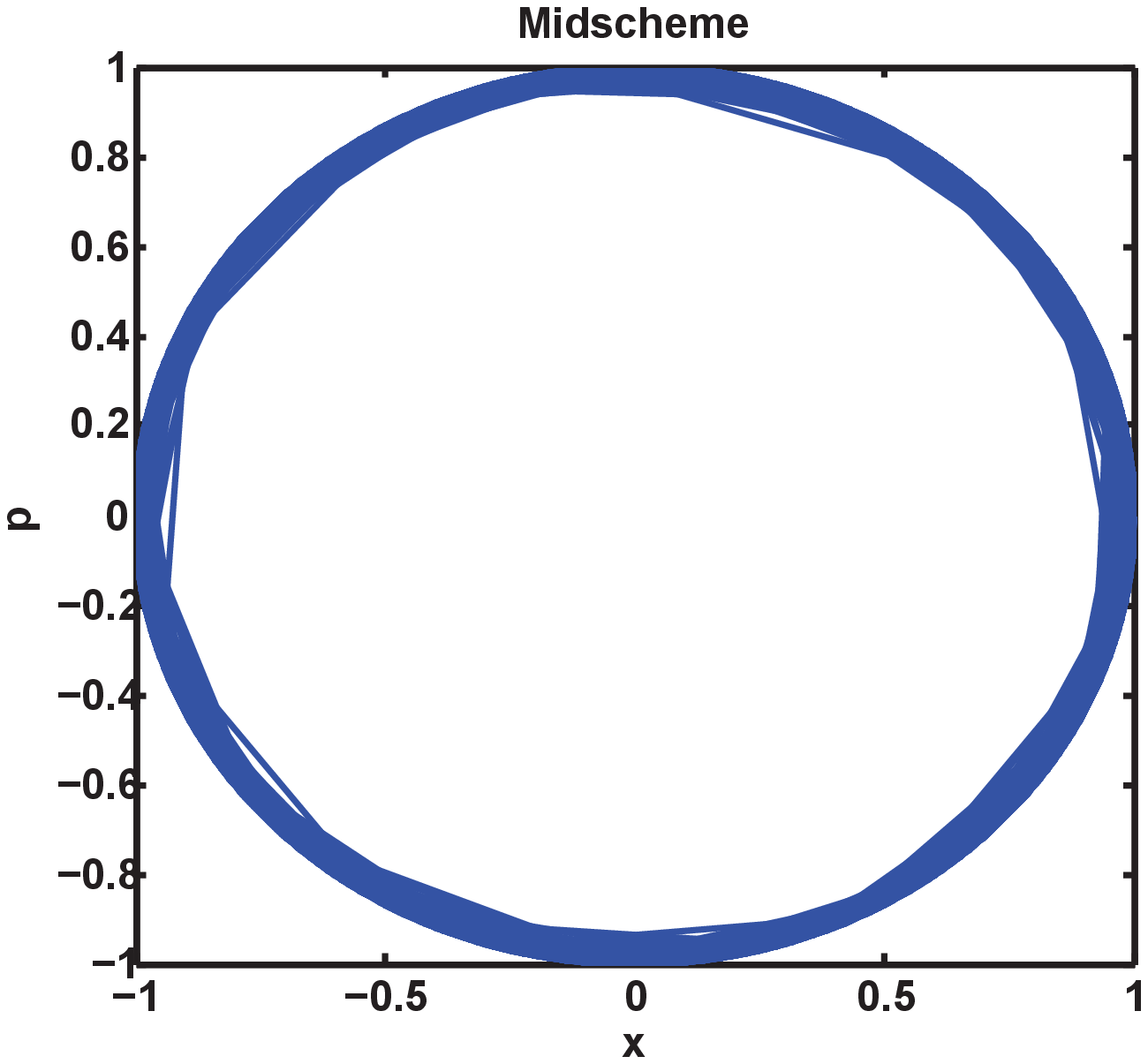}
\includegraphics[height=5cm,width=5cm]{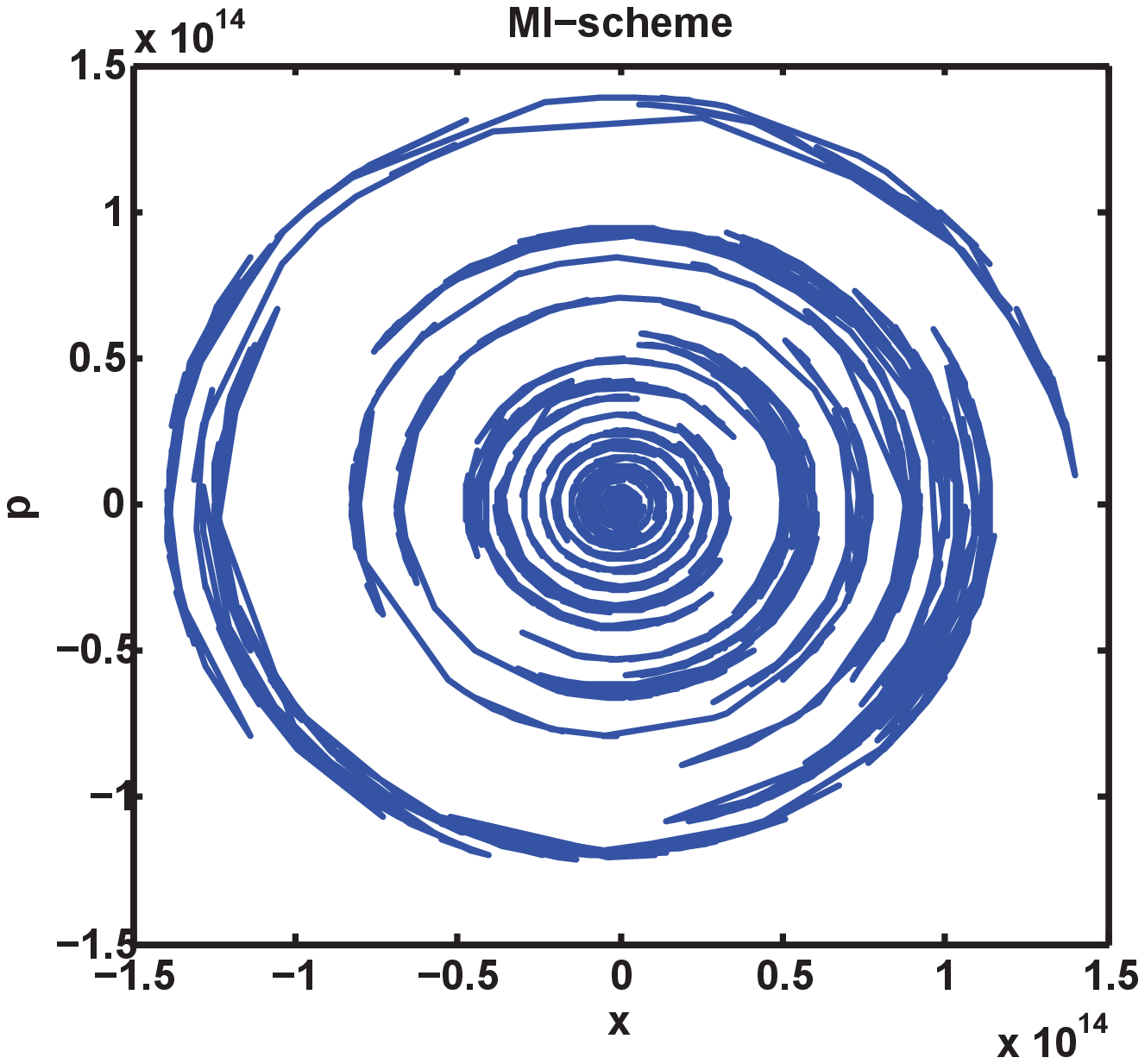}
\caption{Numerical phase trajectory produced by Scheme 2.1 (top left), Scheme 2.2 (top right), midpoint scheme (bottom left) and  Milstein scheme (bottom right) for $\sigma=2$.}
\label{figure05}
\end{center}
\end{figure}

Set $h=0.01$, $T=500$ and $(p_0,x_0)=(0,1)$. We do experiments for $\sigma=1$ and $\sigma=2$.
It can be seen from Figures \ref{figure03} and \ref{figure05} that, after a long period of time, the numerical solution produced by the  Milstein scheme differs gradually form
the initial circle, the numerical solution produced by midpoint scheme runs along the energy circle, and the numerical solutions obtained by Schemes \ref{scheme1} and \ref{scheme2} run
near the initial circle. These phenomenons indicate that Schemes \ref{scheme1} and \ref{scheme2} preserve the quadratic invariant $H(p,x)$ better than the Milstein scheme and worse than
the midpoint scheme. Meanwhile, for small values of $\sigma$ (such as $\sigma\leq 1.0$), Schemes \ref{scheme1} and \ref{scheme2} preform as well as the midpoint scheme in
preserving quadratic invariants, and hence the midpoint scheme can be considered as good templates of implicit conservative methods.

\section{Numerical Tests}
\subsection{Fixed-point iteration}
\label{sec:2}
Despite the actual explicitness of the midpoint rule applied to the Kubo oscillator (\ref{kubo}), we treat it as an usual implicit method with
truncation of the Wiener increments and fixed-point iterations. It takes the form
\begin{equation}\label{mid}\begin{split}
p_{n+1}&=p_n-ah\frac{x_n+x_{n+1}}{2}-\sigma \overline{\Delta W_n}\frac{x_n+x_{n+1}}{2},\quad p_0=1,\\
x_{n+1}&=x_n+ah\frac{p_n+p_{n+1}}{2}+\sigma \overline{\Delta W_n}\frac{p_n+p_{n+1}}{2},\quad x_0=0.
\end{split}
\end{equation}

The fixed-point iteration applied to (\ref{mid}) reads
\begin{equation}\label{midfix}
\begin{split}
p_{n+1}^{[0]}&=p_{n},\quad x_{n+1}^{[0]}=x_n,\quad p_0=1,\quad x_0=0,\\
p_{n+1}^{[N]}&=p_n-ah\frac{x_n+x_{n+1}^{[N-1]}}{2}-\sigma \overline{\Delta W_n}\frac{x_n+x_{n+1}^{[N-1]}}{2},\\
x_{n+1}^{[N]}&=x_n+ah\frac{p_n+p_{n+1}^{[N-1]}}{2}+\sigma \overline{\Delta W_n}\frac{p_n+p_{n+1}^{[N-1]}}{2},\quad N=1,2,\cdots.
\end{split}
\end{equation}

In the numerical tests, we observe the effect on the preservation of the quadratic invariant (\ref{circle}) when taking
different choices of the iteration number $N$, the step-size $h$, and the terminal time $T$ , as well as the convergence rate of the quadratic invariant. In the following, we take $a=2$ and  $\sigma=0.3$. Since the midpoint rule applied to the Kubo oscillator is of root-mean-square order 1 (see \cite{3}), we take $k=2$ in the truncation $\overline{\Delta W_n}$ of realizing $\Delta W_n$.

\begin{figure}[htp]
\begin{center}
    \includegraphics[width=5cm,height=5cm]{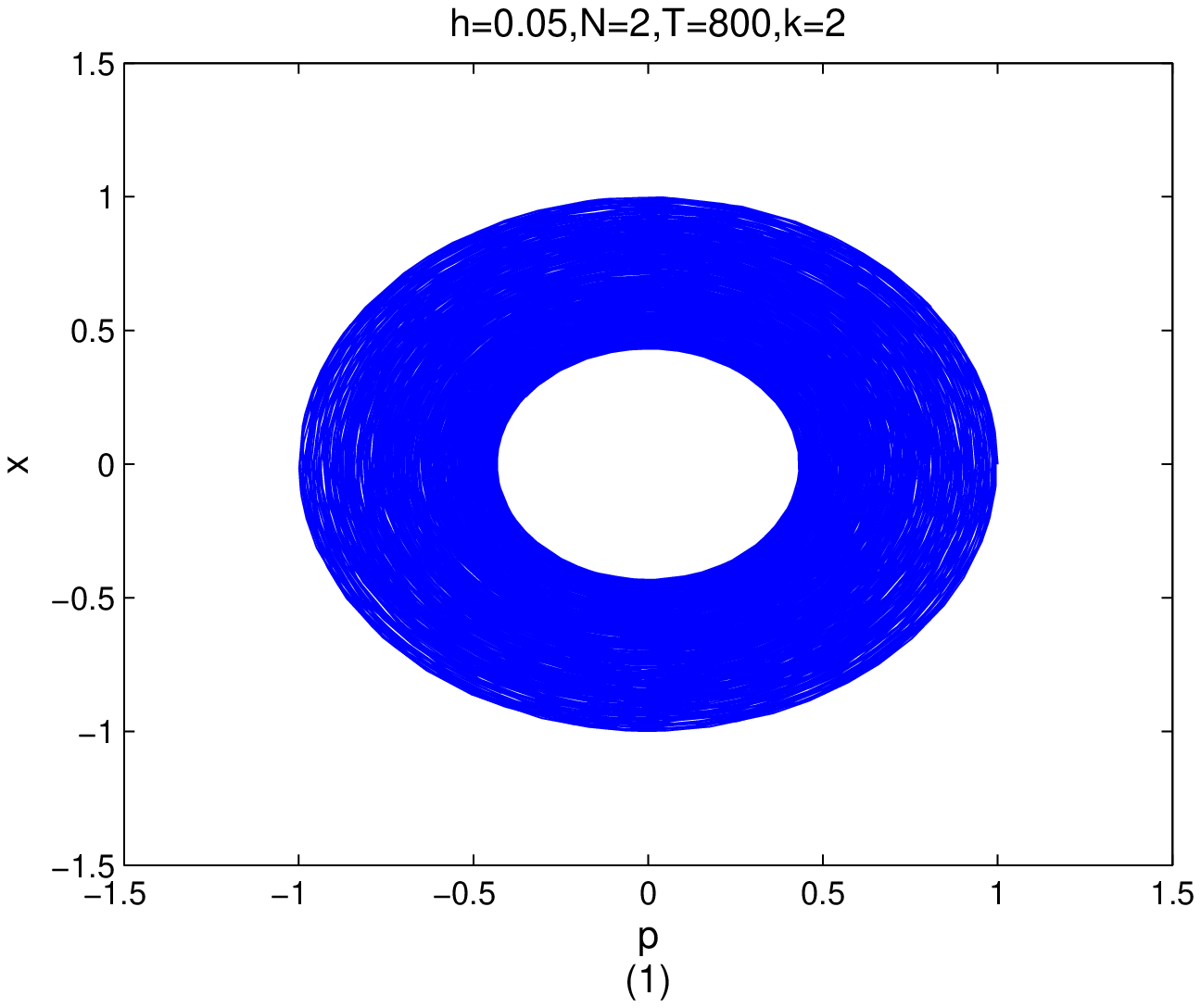}
    \includegraphics[width=5cm,height=5cm]{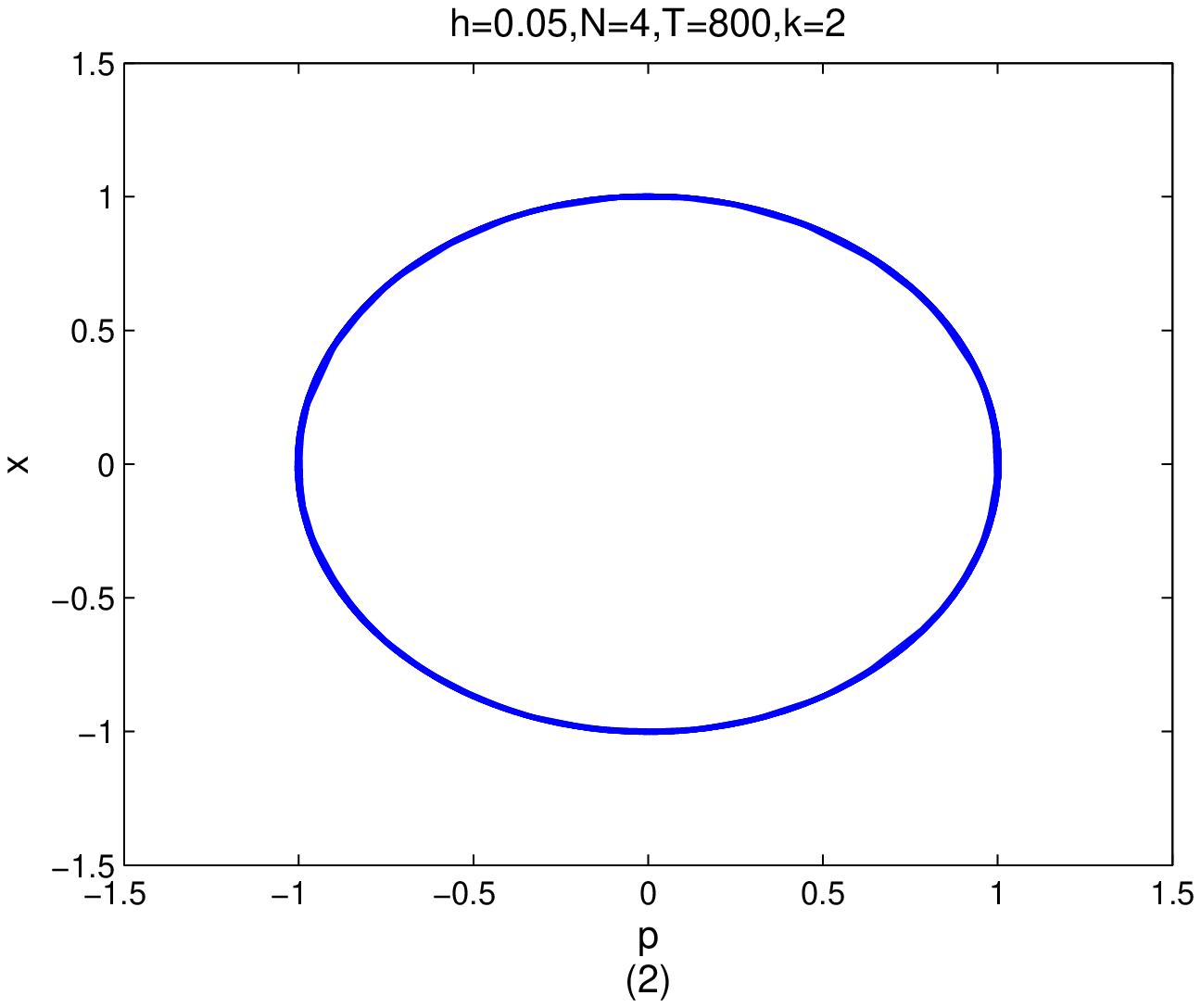}
\caption{Phase trajectory produced by (\ref{midfix}) for (1) $N=2$ (2) $N=4$, $h=0.05$ and $t\in[0,800]$.}\label{f1}
  \end{center}
\end{figure}
\begin{figure}[htp]
\begin{center}
       \includegraphics[width=5cm,height=5cm]{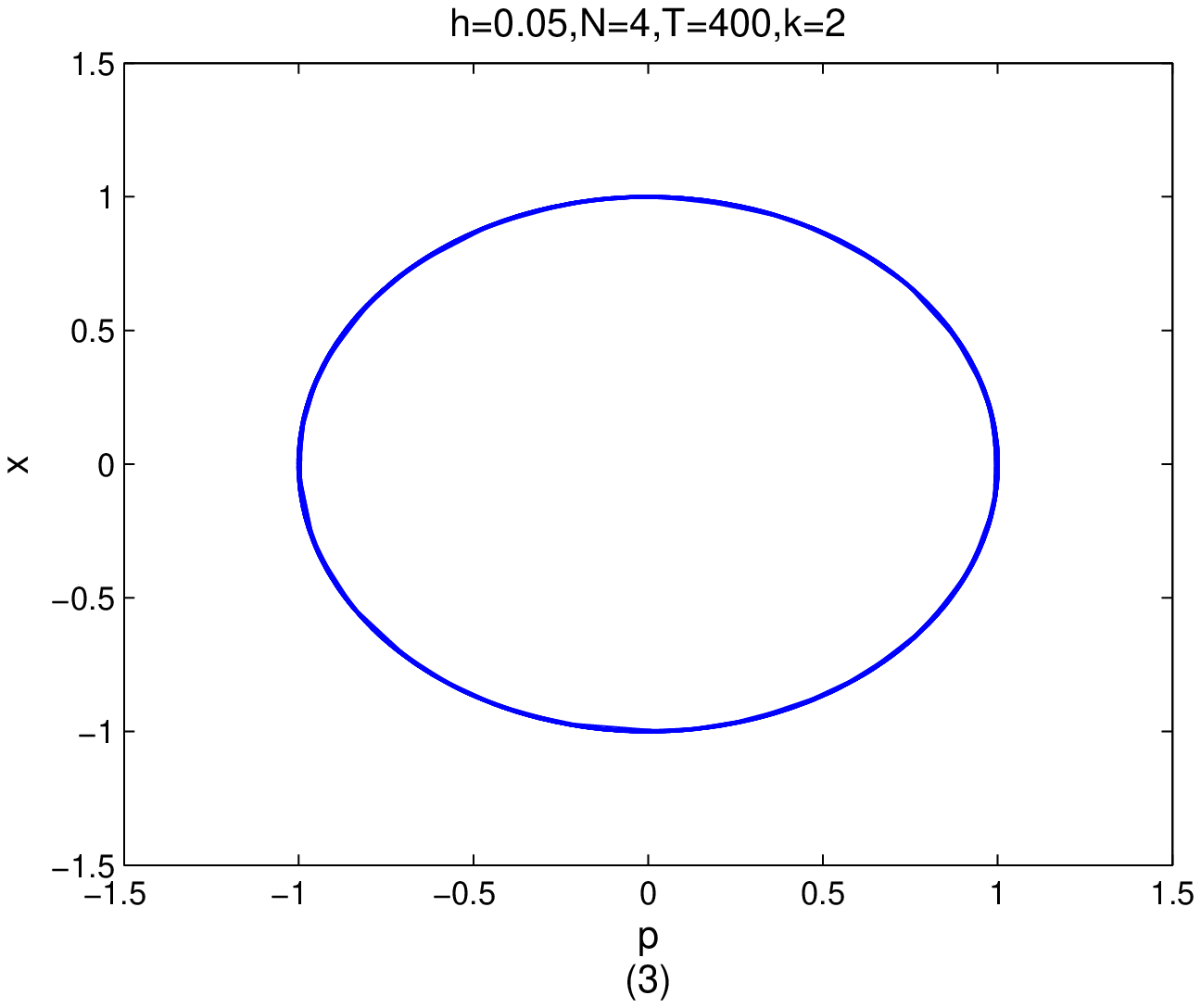}
       \includegraphics[width=5cm,height=5cm]{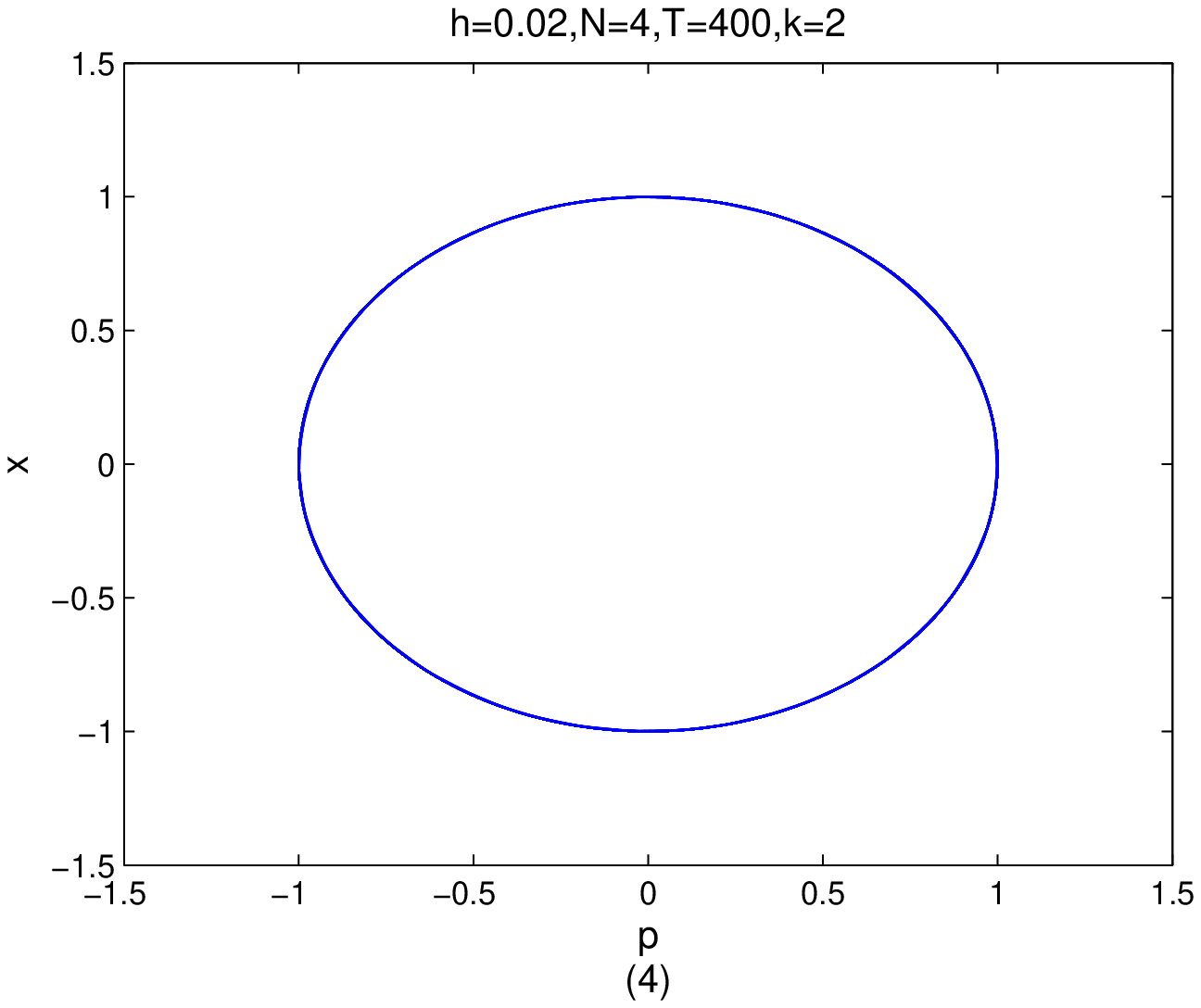}
\caption{Phase trajectory produced by (\ref{midfix}) for (3) $h=0.05$ (4) $h=0.02$, $N=4$ and $t\in[0,400]$.}\label{f2}
\end{center}
\end{figure}
\begin{figure}[htp]
\begin{center}
       \includegraphics[width=5cm,height=5cm]{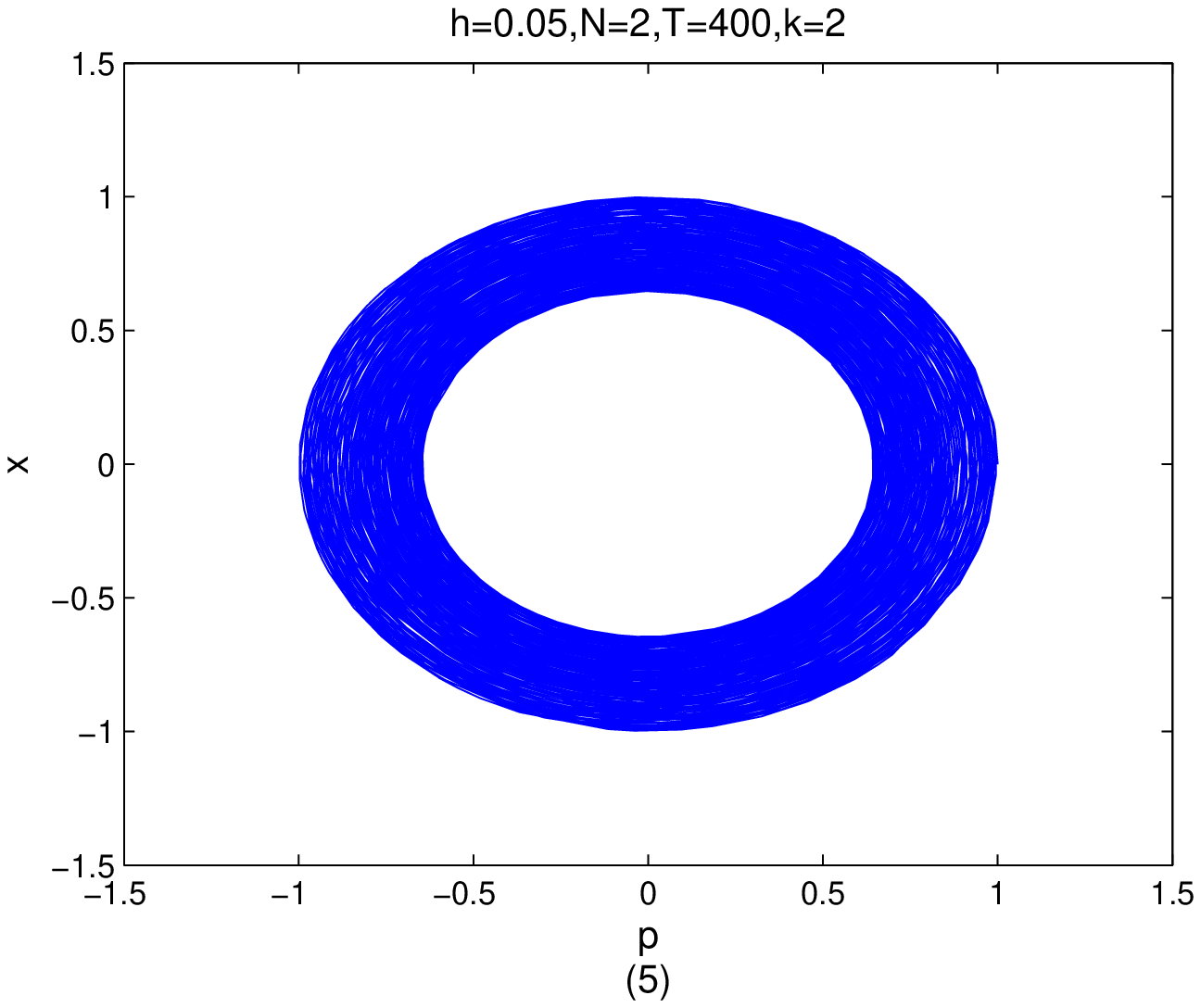}
       \includegraphics[width=5cm,height=5cm]{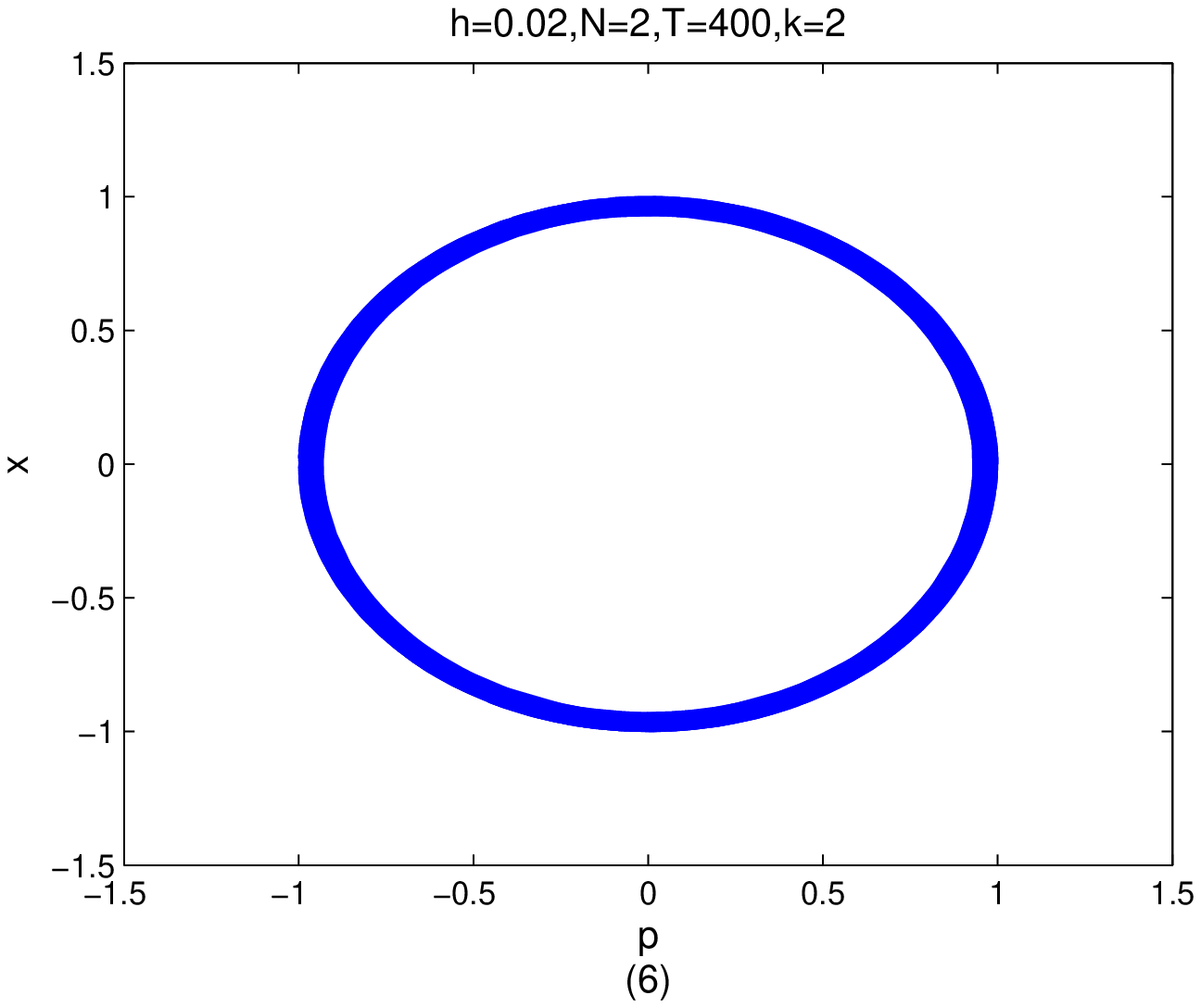}
\caption{Phase trajectory produced by (\ref{midfix}) for (5) $h=0.05$ (6) $h=0.02$, $N=2$ and $t\in[0,400]$. }\label{f3}
\end{center}
\end{figure}
\begin{figure}[htp]
\begin{center}
    \includegraphics[width=5cm,height=5cm]{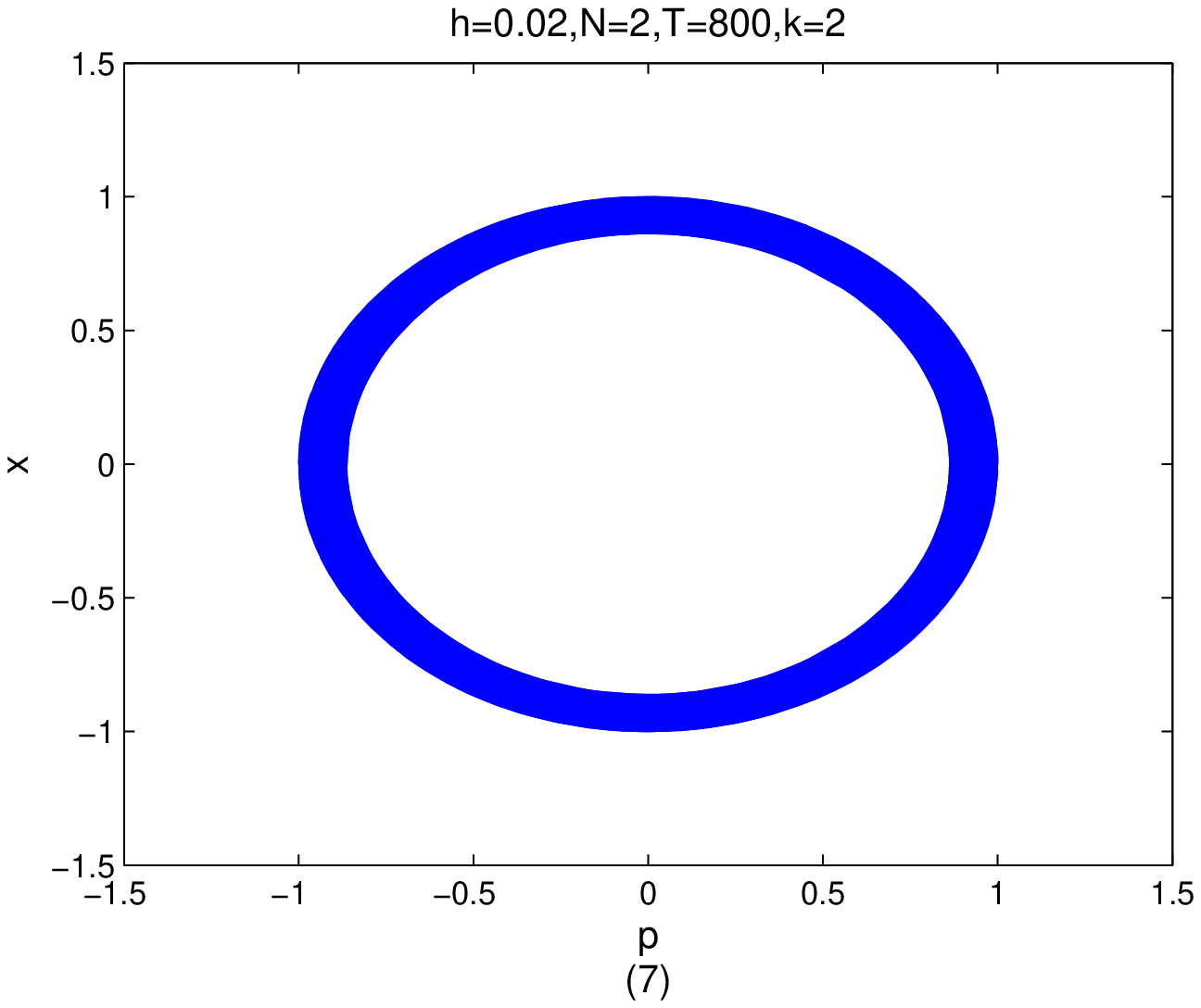}
    \includegraphics[width=5cm,height=5cm]{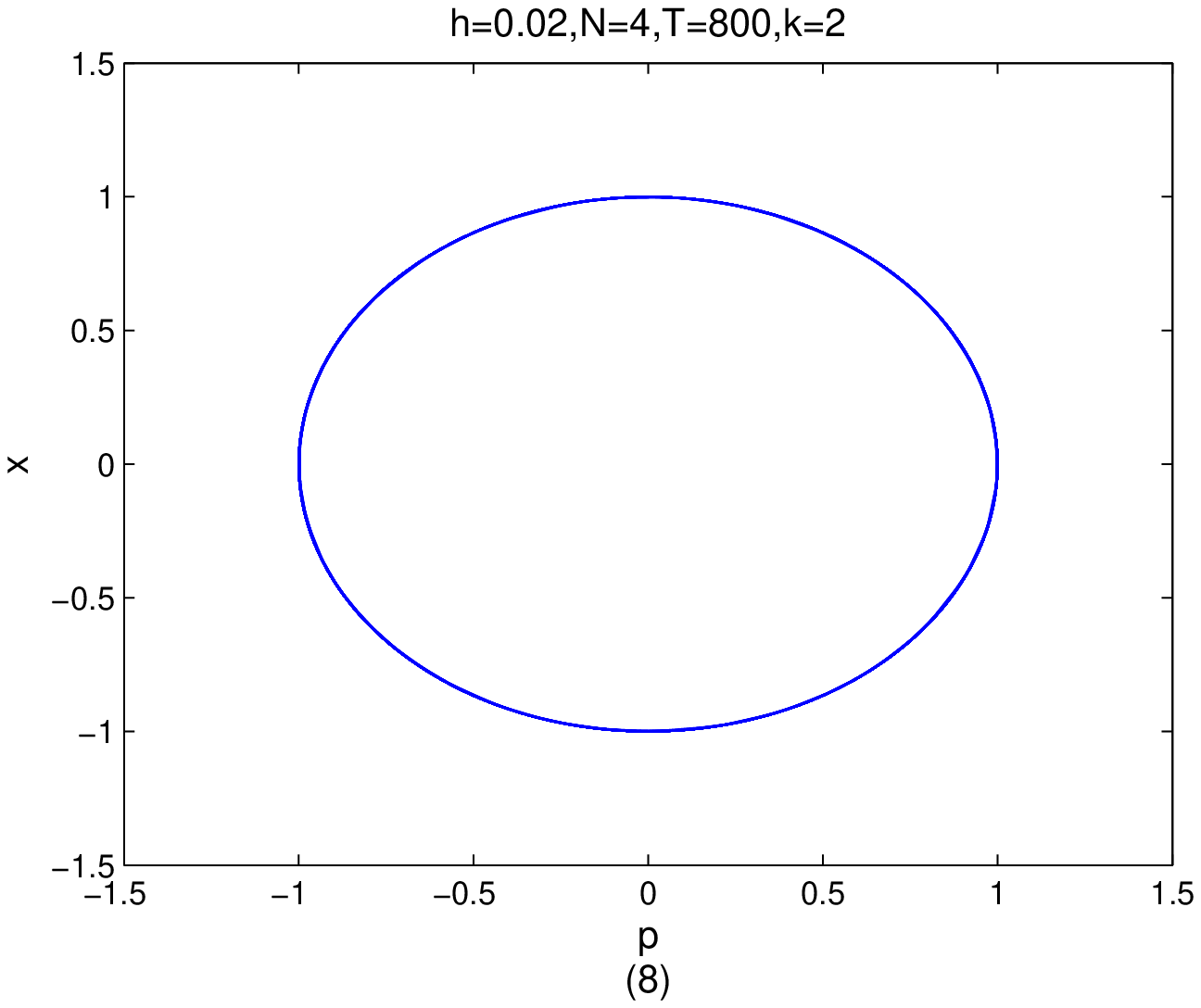}
\caption{Phase trajectory produced by (\ref{midfix}) for (7) $N=2$ (8) $N=4$, $h=0.02$ and $t\in[0,800]$.}\label{f4}
\end{center}
\end{figure}

Figures \ref{f1} to \ref{f4} compare the phase trajectories produced by (\ref{midfix}) with different choices of iteration number $N$, time step-size $h$, and terminal time $T$. It can be seen from Figures \ref{f1} and \ref{f4} that, as $N$ increases from $2$ to $4$, the quadratic invariant is much better preserved, while for a smaller $h$ ($h=0.02$), the effect is even better. Under the same setting of $N$ and $T$, the preservation of quadratic invariant is more accurate for $h=0.02$ than $h=0.05$, as can be observed from Figures \ref{f2} and \ref{f3}, by which it is also indicated that $N=4$ gives better preservation of quadratic invariant than $N=2$. As far as $T$ is concerned, the panels (1) and (5) ($h=0.05$) show the accumulation of the error in the quadratic invariants as $T$ gets larger, to which the panels (6) and (7) for $h=0.02$ are contributed as well.
\begin{figure}[htp]
\begin{center}
    \includegraphics[width=3cm,height=5cm]{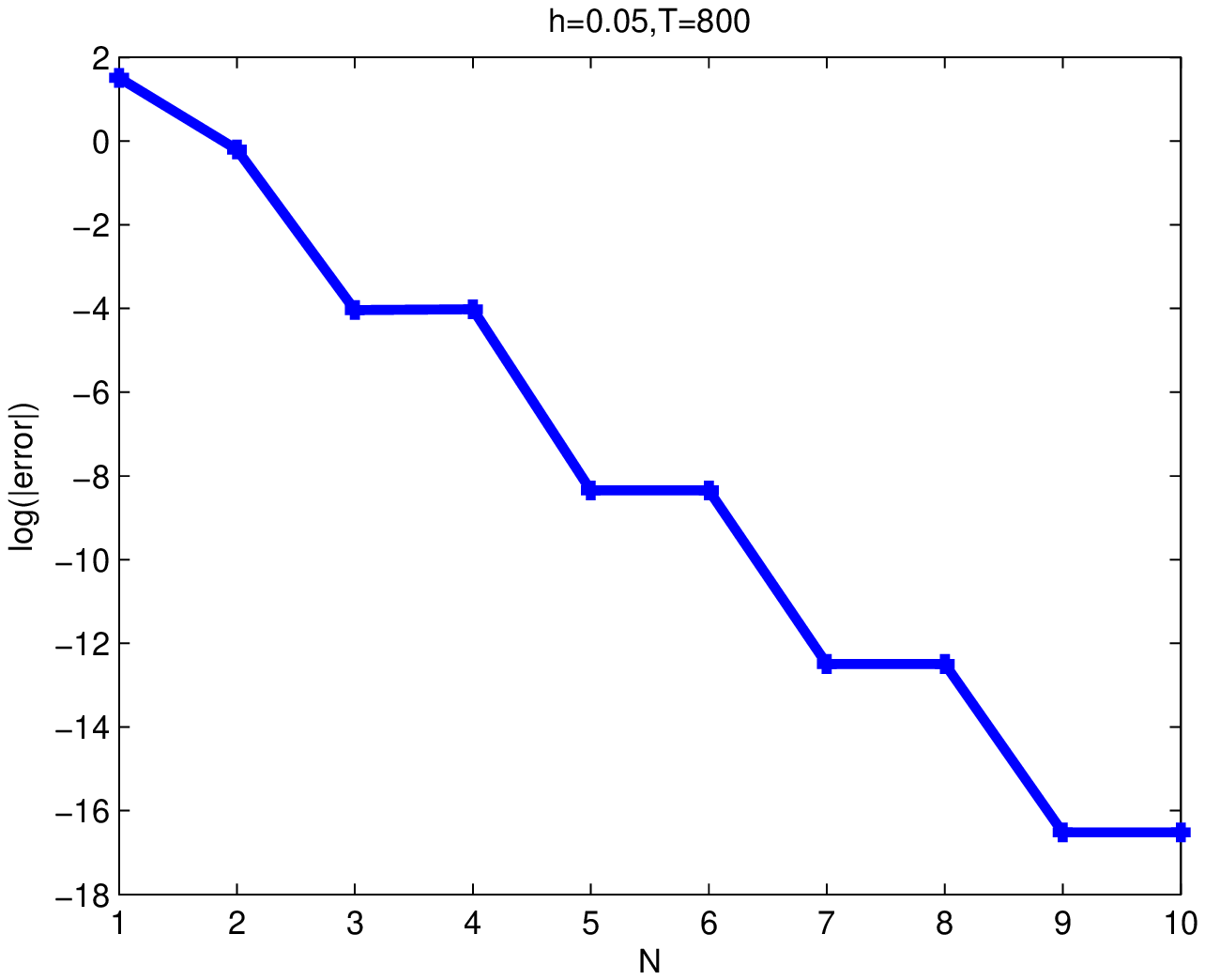}
    \includegraphics[width=3cm,height=5cm]{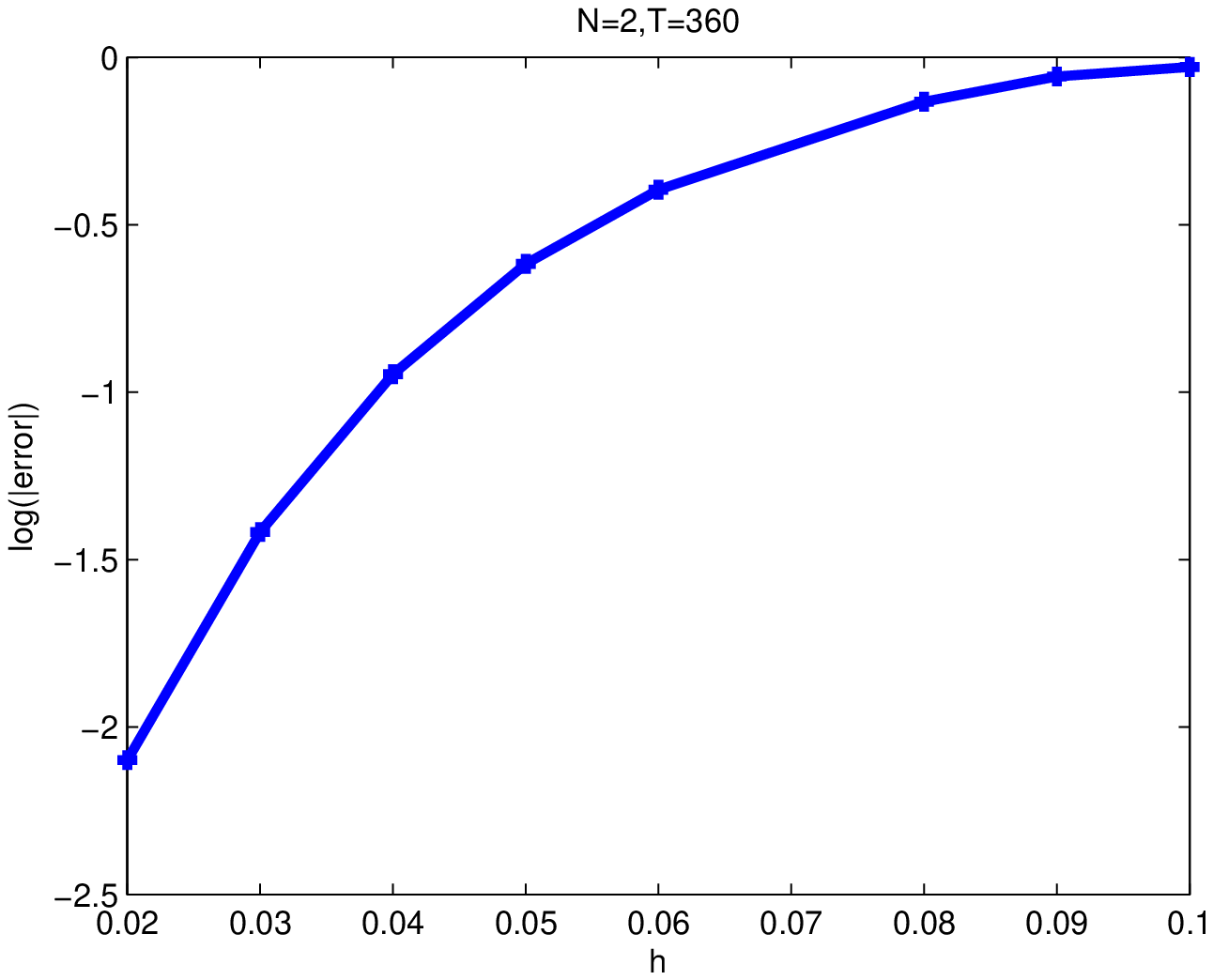}
    \includegraphics[width=3cm,height=5cm]{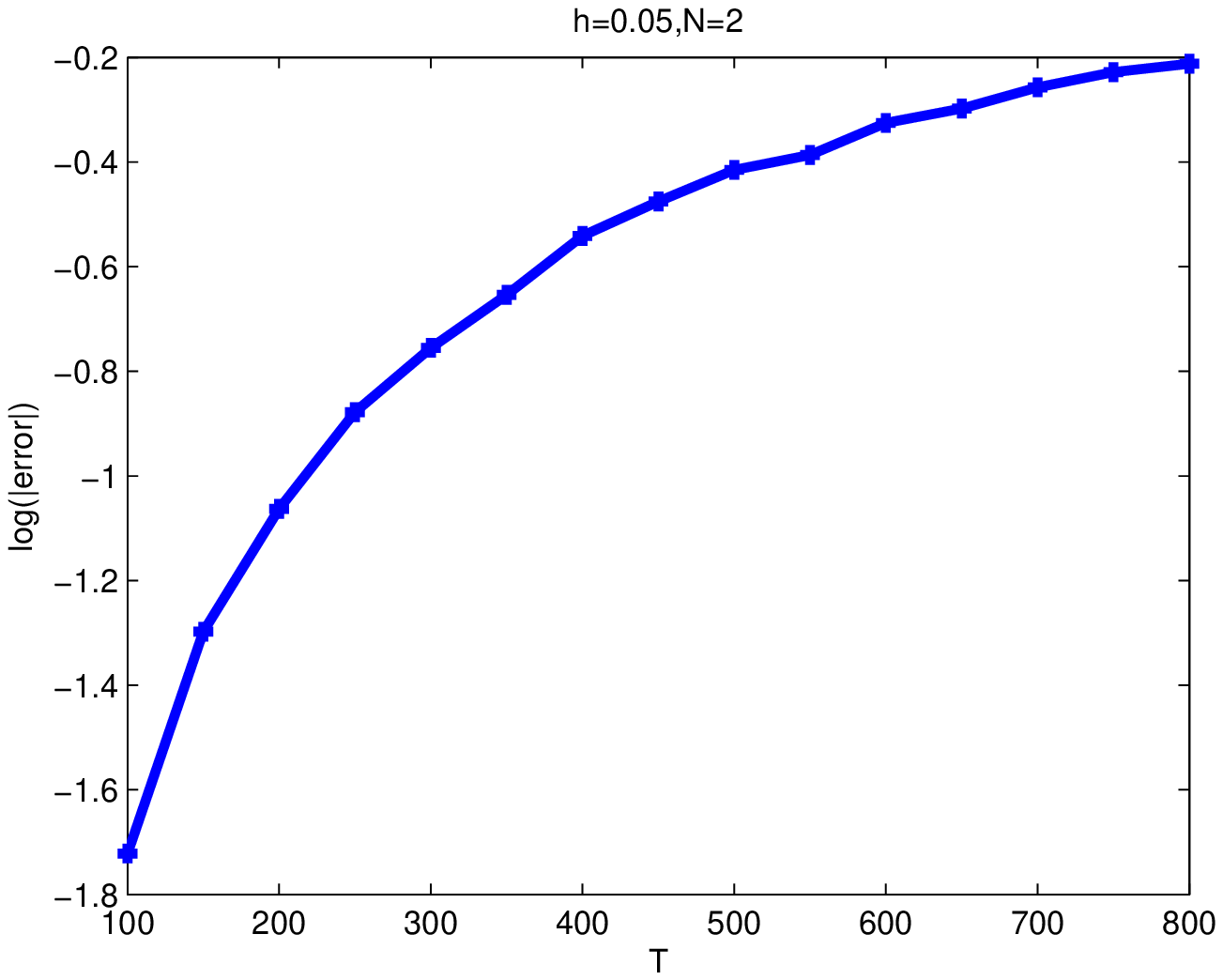}
\caption{Logarithm of the error in the quadratic invariant against (1) $N$, for $N\in\{1,2,3,\cdots, 10\}$, $h=0.05$, $T=800$ (left); (2) $h$, for $h\in\{0.02, 0.03, 0.04, 0.05, 0.06, 0.08, 0.09, 0.1\}$, $N=2$, $T=360$ (middle); (3) $T$, for $T\in\{100,150,200,\cdots,800\}$, $h=0.05$, $N=2$ (right).}\label{f5}
  \end{center}
\end{figure}

Suppose $V=\frac{T}{h}$. Figure \ref{f5} is devoted to illustrating the change of logarithm of the error in the quadratic invariant at time $T$ $$\ln|(p_{V}^{[N]})^2+(q_V^{[N]})^2-1|$$ against that of $N$, $h$ and $T$, from which the decrease of the error with the increase of $N$, as well as the increase of the error with the increase of $h$ and $T$ can be easily seen.

\begin{figure}[htp]
\begin{center}
    \includegraphics[width=5cm,height=5cm]{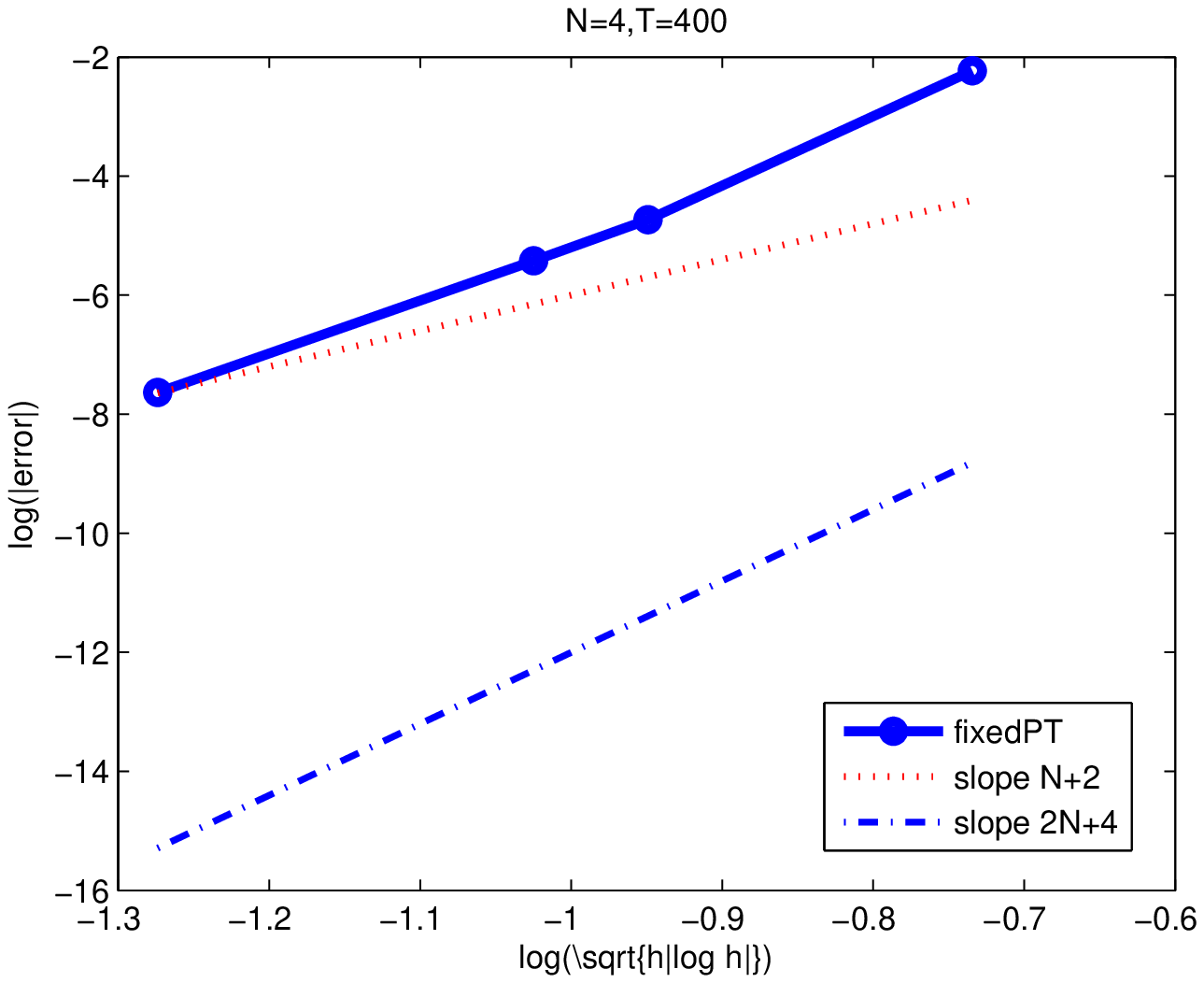}
    \includegraphics[width=5cm,height=5cm]{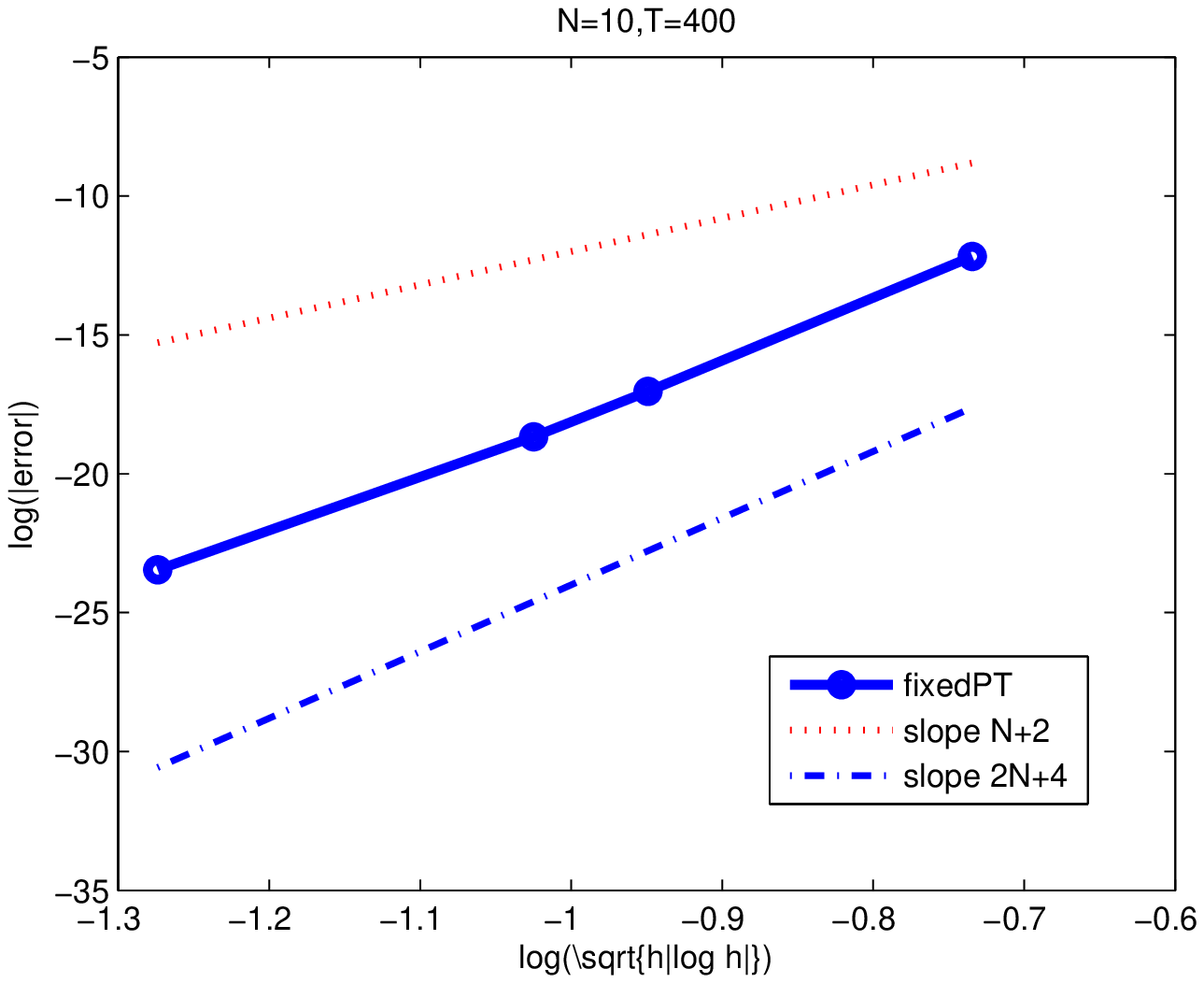}
\caption{Logarithm of error in the quadratic invariant against logarithm of $\sqrt{h|\ln h|}$ for $h\in\{0.01,0.02,0.04,0.08\}$, $T=400$, and (1) $N=4$ (left); (2) $N=10$ (right).}\label{f6}
  \end{center}
\end{figure}
\begin{figure}[htp]
\begin{center}
    \includegraphics[width=5cm,height=5cm]{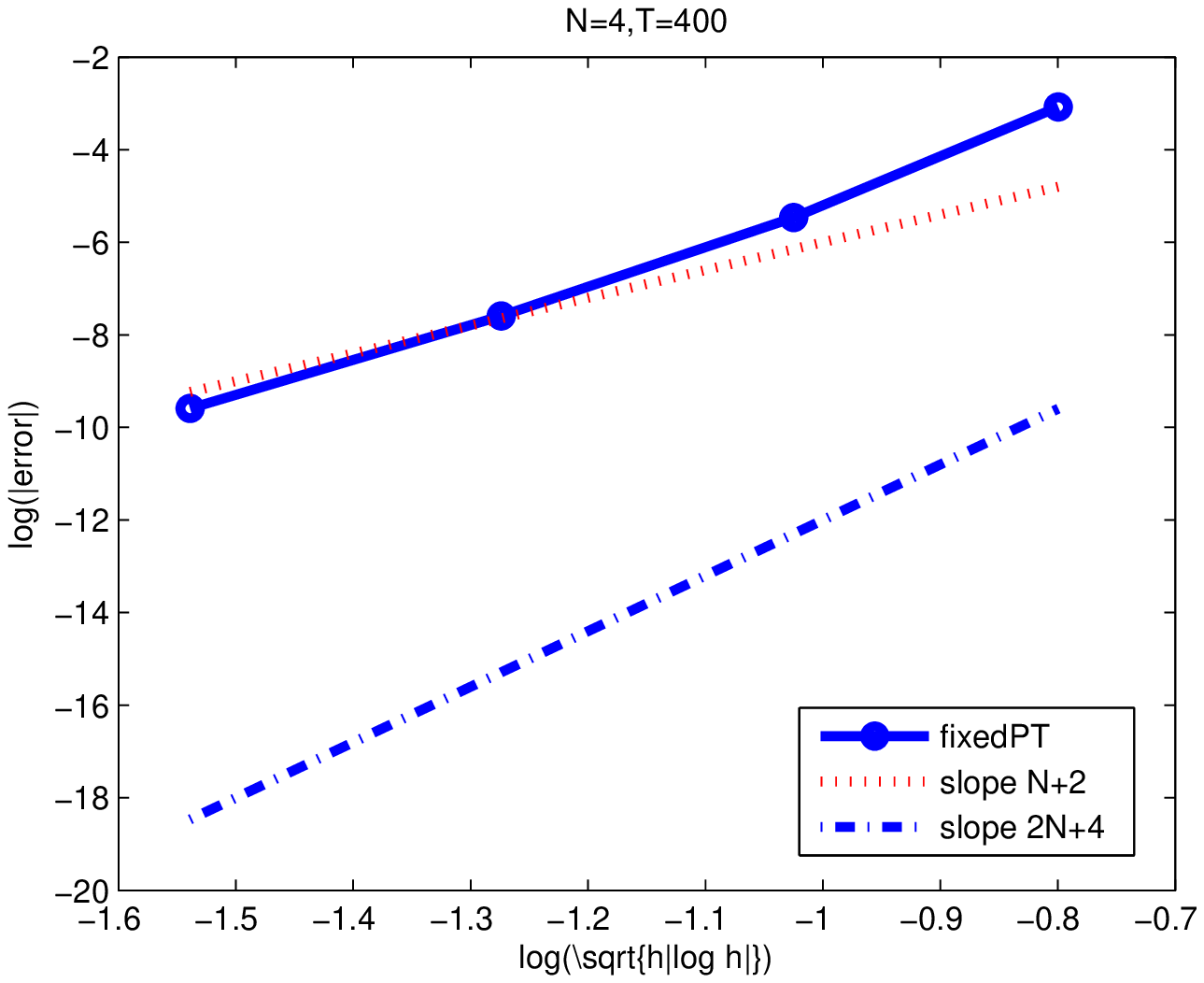}
    \includegraphics[width=5cm,height=5cm]{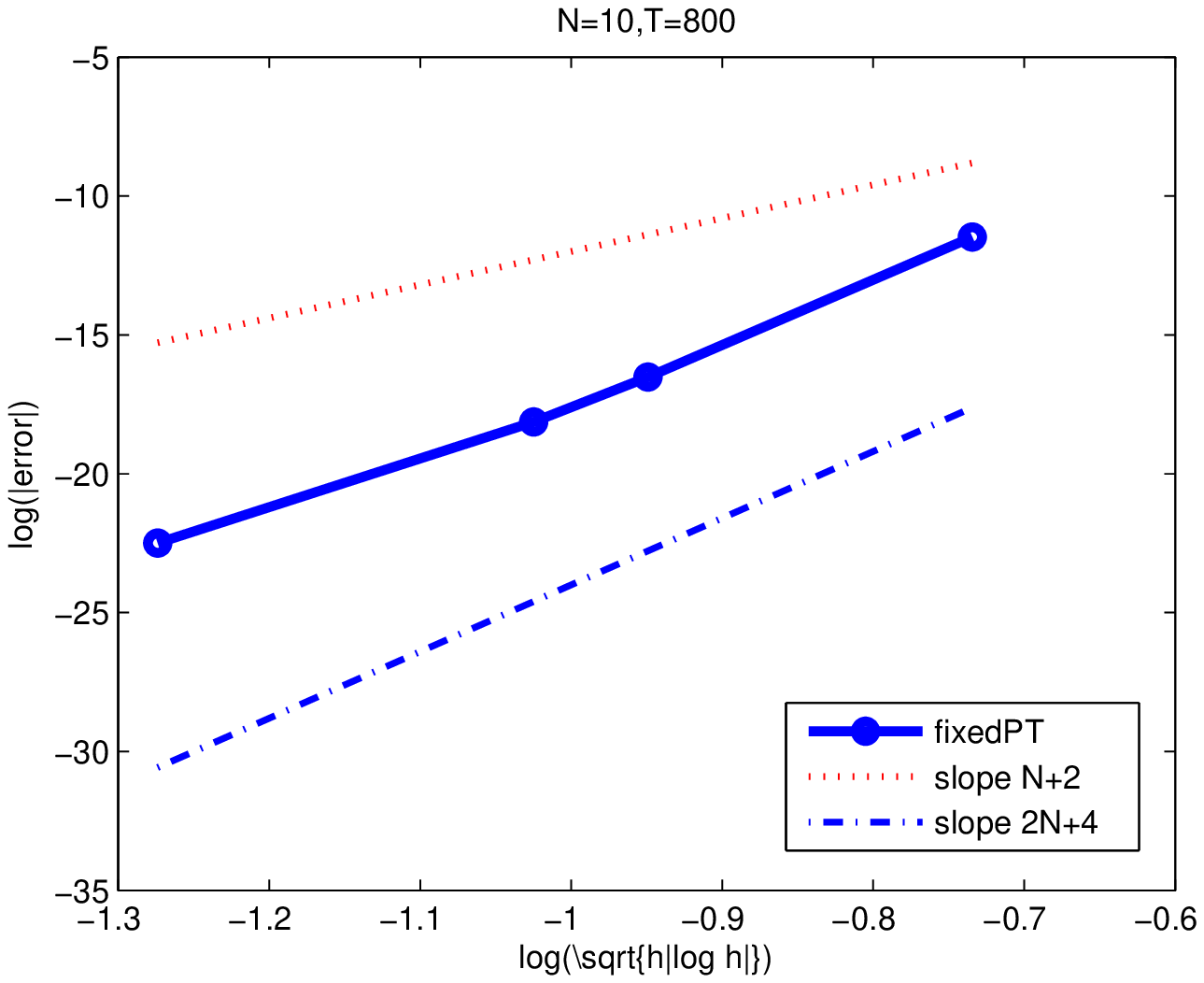}
\caption{Logarithm of error in the quadratic invariant against logarithm of $\sqrt{h|\ln h|}$ for $h\in\{0.01,0.02,0.04,0.08\}$, $T=800$, and (1) $N=4$ (left); (2) $N=10$ (right). }\label{f7}
  \end{center}
\end{figure}

The convergence rate of the error in the preservation of quadratic invariant by the fixed-point iteration, which is theoretically indicated in Theorem \ref{fixed thm2}, is numerically demonstrated via Figures \ref{f6} and \ref{f7}, where the red dotted line is a reference of slope $N+2$, and the blue dash-dotted line is a reference of slope $2N+4$ in both figures. The overall speed of convergence seems to be between $O(\sqrt{h|\ln h|}^{N+2})$ and $O(\sqrt{h|\ln h|}^{2N+4})$, as given in Theorem \ref{fixed thm2}. Meanwhile, as $T$ gets larger, the slope of the numerical line is invariant, but with a slight parallel upward translation.

\subsection{Newton's iteration}
\label{sec:2}
Due to the linearity of the system (\ref{kubo}), the Newton's iteration applied to the midpoint rule (\ref{mid}) for the Kubo oscillator (\ref{kubo}) reverts to the method (\ref{mid}) itself.
So we need a nonlinear system to test the behavior of the Newton's iteration. Consider the following stochastic Hamiltonian system
\begin{equation}\label{testSHS}
\begin{split}
dp&=-pqdt-\frac{1}{2}p^2\circ dW(t),\quad p(0)=p_0,\\
dq&=\frac{1}{2}q^2dt+pq\circ dW(t),\quad q(0)=q_0,
\end{split}
\end{equation}
with Hamiltonians $H_0=\frac{1}{2}pq^2$ and $H_{1}=\frac{1}{2}p^2q$. It is easy to verify that $dp(t)\wedge dq(t)=dp_0\wedge dq_0$, $\forall t\geq 0$, which means geometrically `area preservation' in the phase space.

The midpoint rule, which is a symplectic method, applied to (\ref{testSHS}) is
\begin{equation}\label{midshs}
\begin{split}
p_{n+1}&=p_n-h(\frac{p_n+p_{n+1}}{2})(\frac{q_n+q_{n+1}}{2})-\frac{1}{2}\overline{\Delta W_n}(\frac{p_n+p_{n+1}}{2})^2,\\
q_{n+1}&=q_n+\frac{h}{2}(\frac{q_n+q_{n+1}}{2})^2+\overline{\Delta W_n}(\frac{p_n+p_{n+1}}{2})(\frac{q_n+q_{n+1}}{2}).
\end{split}
\end{equation}
The Newton's iteration for solving $(p_{n+1},q_{n+1})$ from the implicit scheme (\ref{midshs}) reads
\begin{equation}\label{newtshs}
\begin{split}
p_{n+1}^{[N]}&=p_{n+1}^{[N-1]}-(\frac{a_n\alpha_n+c_n\beta_n}{a_nb_n-c_nd_n})^{[N-1]},\quad p_{n+1}^{[0]}=p_n,\\
p_{n+1}^{[N]}&=q_{n+1}^{[N-1]}-(\frac{d_n\alpha_n+b_n\beta_n}{a_nb_n-c_nd_n})^{[N-1]},\quad q_{n+1}^{[0]}=q_n,\quad N=1,2,\cdots,
\end{split}
\end{equation}
where
\begin{equation}\label{parameters}
\begin{split}
a_n&=1-\frac{hv_n}{4}+\frac{u_n}{4}\overline{\Delta W_n},\quad c_n=-\frac{hu_n}{4},\\
b_n&=1+\frac{hv_n}{4}+\frac{u_n}{4}\overline{\Delta W_n},\quad d_n=-\frac{v_n}{4}\overline{\Delta W_n},\\
\alpha_n&=p_{n+1}-p_n+\frac{h}{4}u_nv_n+\frac{u_n^2}{8}\overline{\Delta W_n},\\
\beta_n&=q_{n+1}-q_n-\frac{h}{8}v_n^2+\frac{u_nv_n}{4}\overline{\Delta W_n},
\end{split}
\end{equation}
with
\begin{equation}\label{alphabeta}
\begin{split}
u_n=p_n+p_{n+1}, \quad v_n=q_{n}+q_{n+1},
\end{split}
\end{equation}
and the upper index $[N-1]$ on the fractions appearing in the right hand side of (\ref{newtshs}) means that all the $p_{n+1}$ and $q_{n+1}$ included in the fractions are with upper index $[N-1]$.

In order to test the area preservation in the phase space, we choose the initial phase points $(p_0,q_0)$ from the unit circle in the $p-q$ plane.
And we observe the evolution of the circle driven by the midpoint rule (\ref{midshs}) with Newton's iteration (\ref{newtshs}) under influence of different choices of iteration number $N$,
time step-size $h$, and terminal time $T$.
\begin{figure}[htp]
\begin{center}
    \includegraphics[width=3cm,height=5cm]{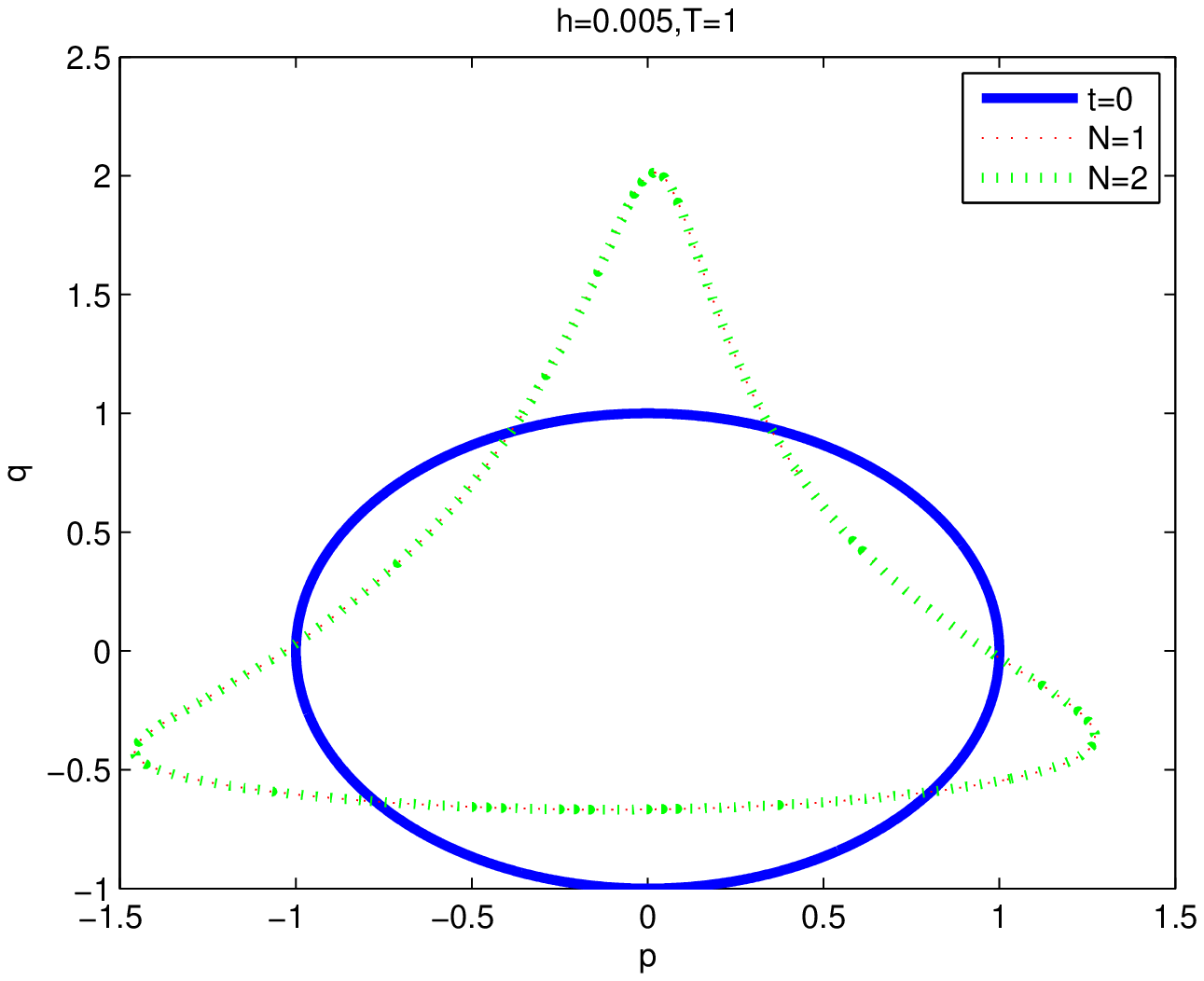}
    \includegraphics[width=3cm,height=5cm]{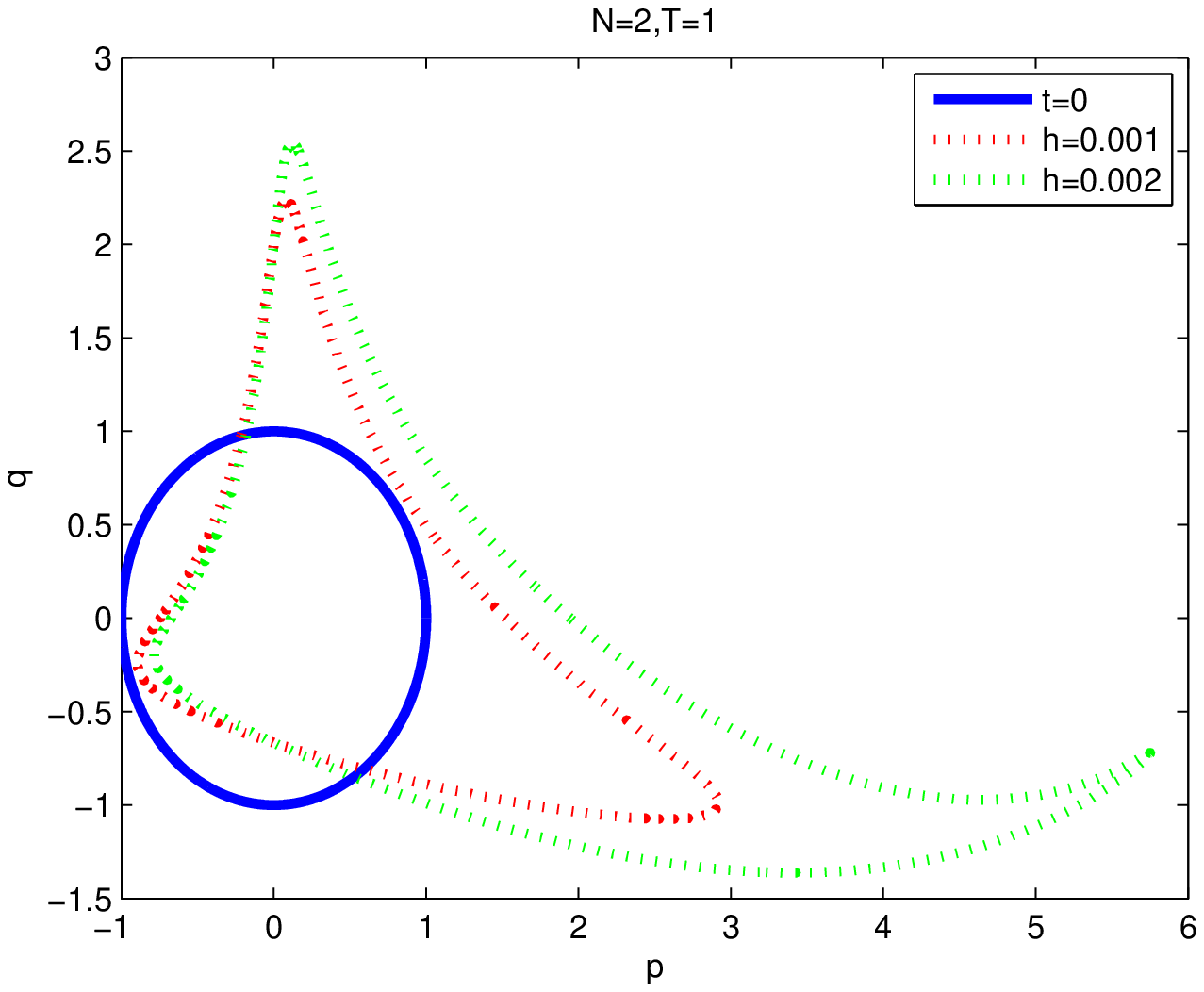}
    \includegraphics[width=3cm,height=5cm]{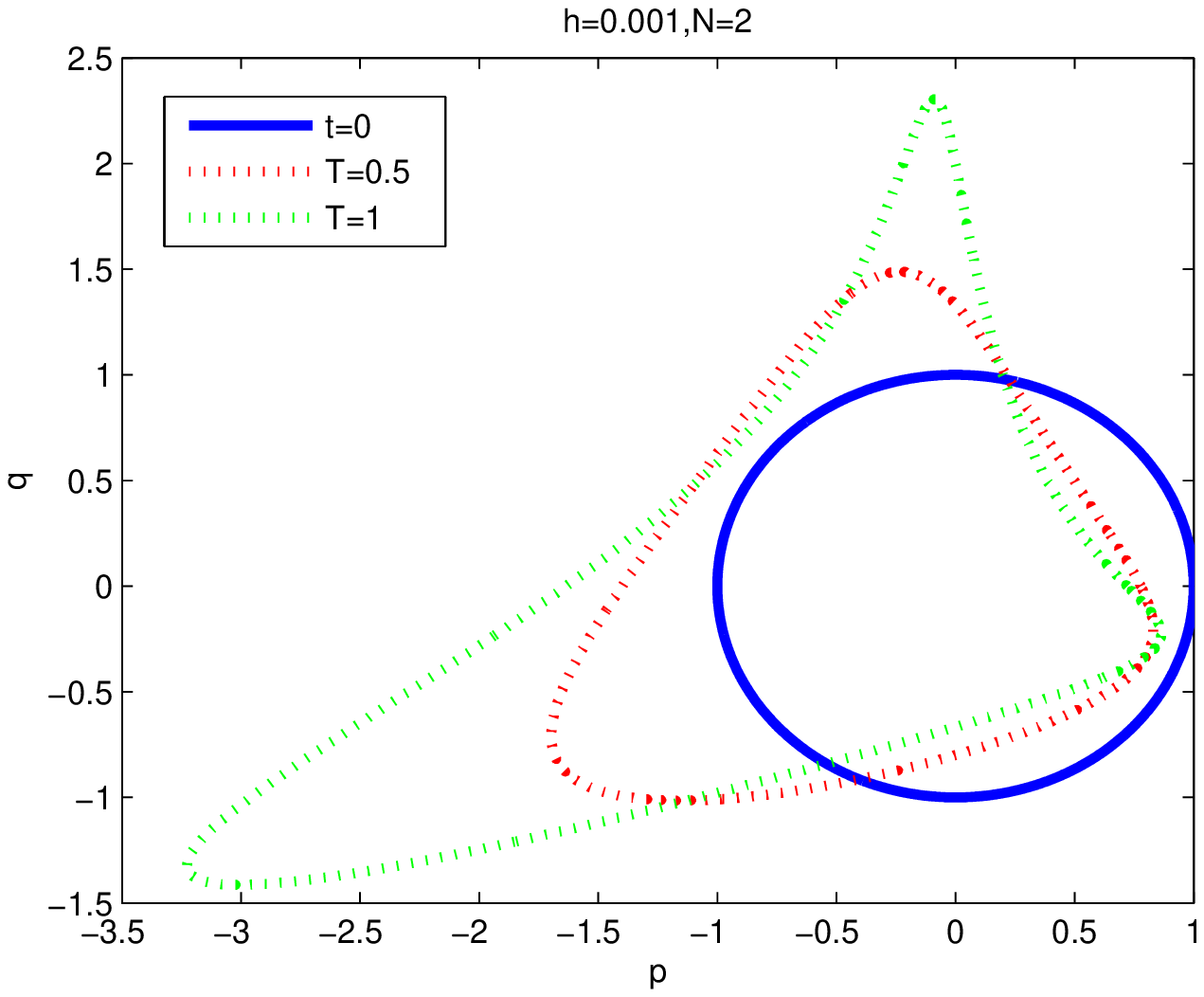}
\caption{Evolution of the unit circle driven by the midpoint rule (\ref{midshs}) with Newton's iteration (\ref{newtshs}) for (left)  $h=0.005$, $T=1$, $N=1$ (red dotted), $N=2$ (green dotted); (middle) $N=2$, $T=1$, $h=0.001$ (red dotted), $h=0.002$ (green dotted); (right) $N=2$, $h=0.001$, $T=0.5$ (red dotted), $T=1$ (green dotted).}\label{f8}
  \end{center}
\end{figure}

It can be seen from the left panel of Figure \ref{f8} that, the red and green dotted lines coincide visually, which indicates that mild variation of the iteration number $N$ within a reasonable domain could have little influence on the area preservation of the algorithm, while the middle and right panels of Figure \ref{f8} illustrate obvious growth of error in the area preservation of the algorithms with the increase of $h$ and $T$, respectively.

\begin{figure}[htp]
\begin{center}
    \includegraphics[width=3cm,height=5cm]{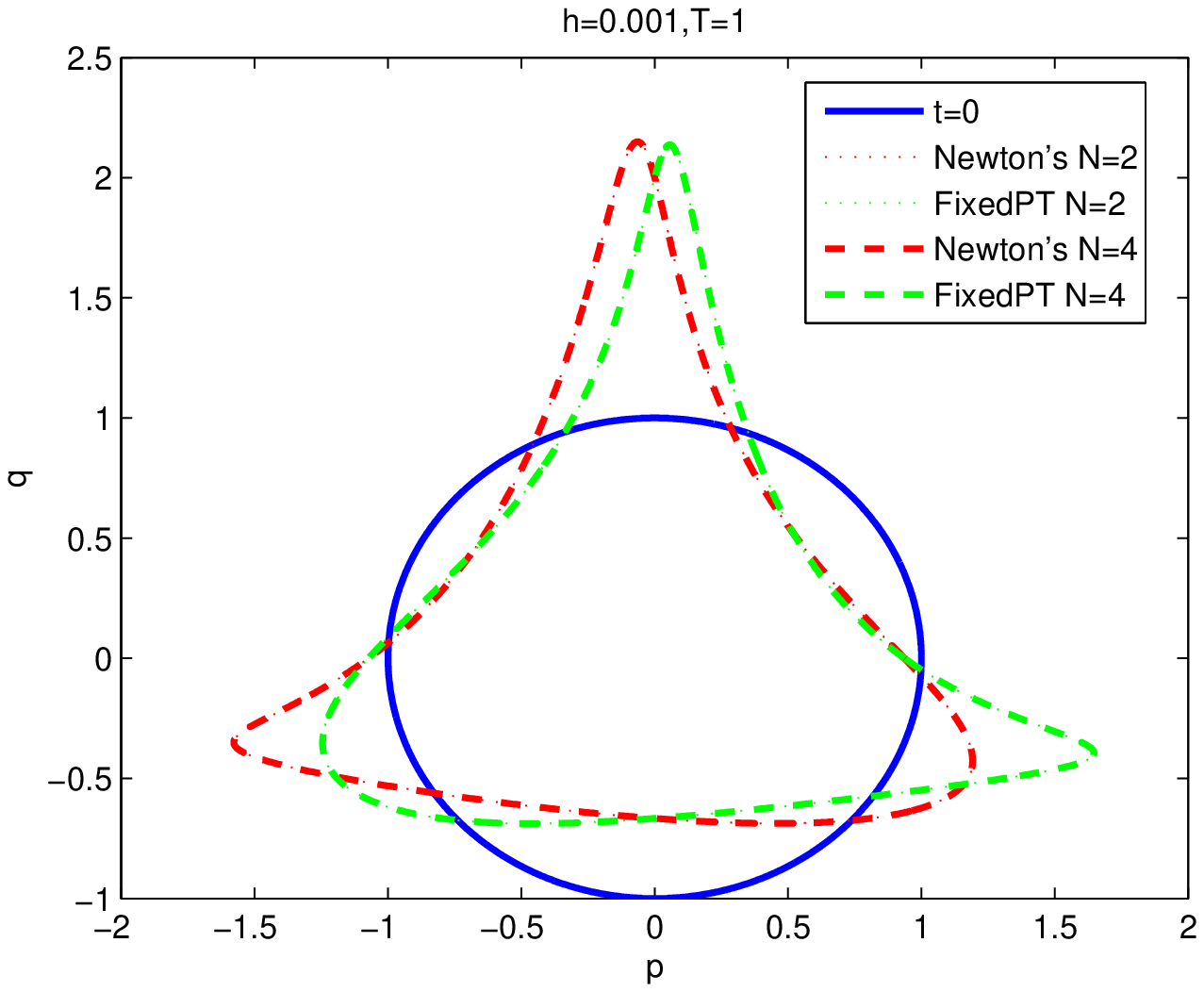}
    \includegraphics[width=3cm,height=5cm]{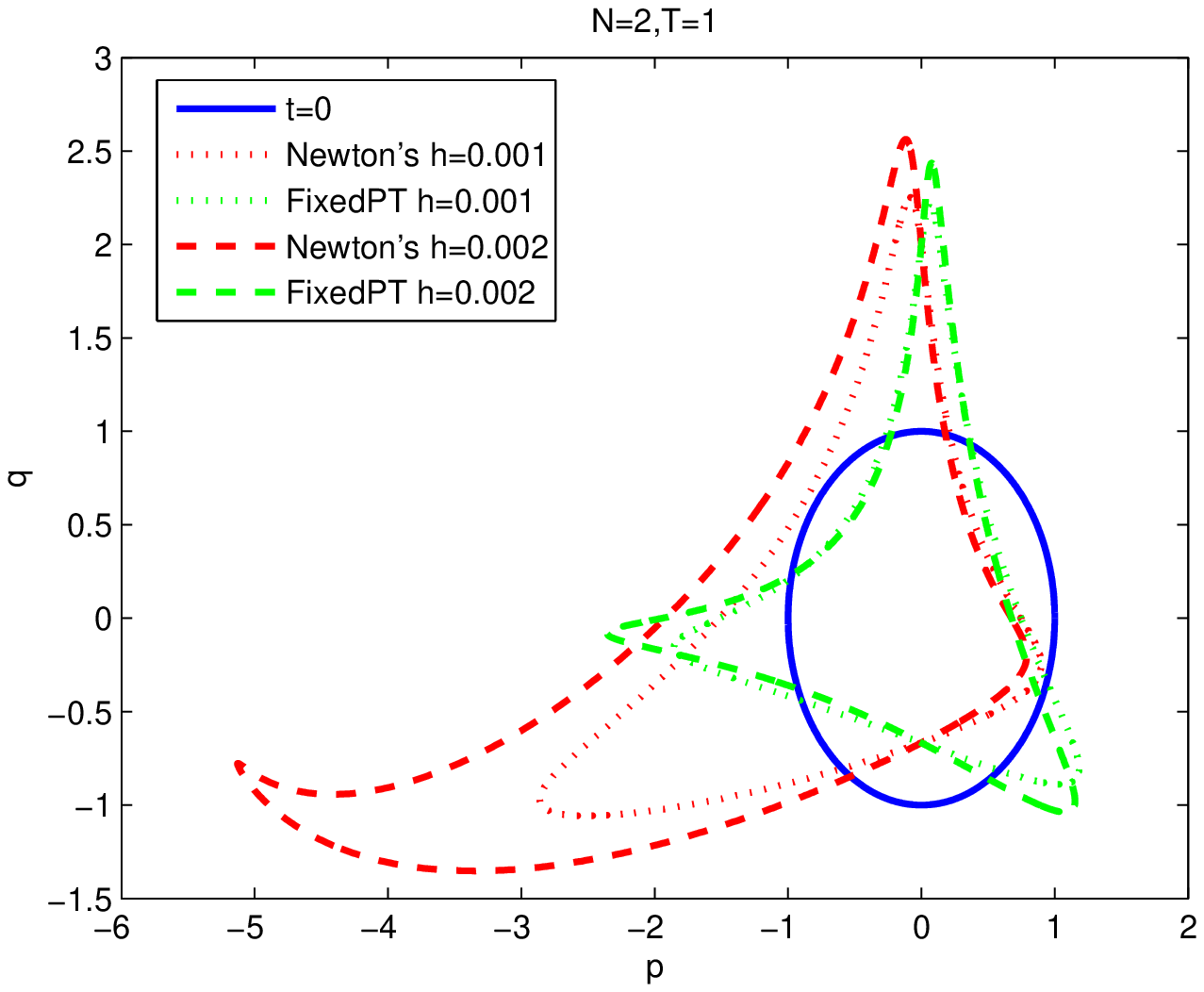}
    \includegraphics[width=3cm,height=5cm]{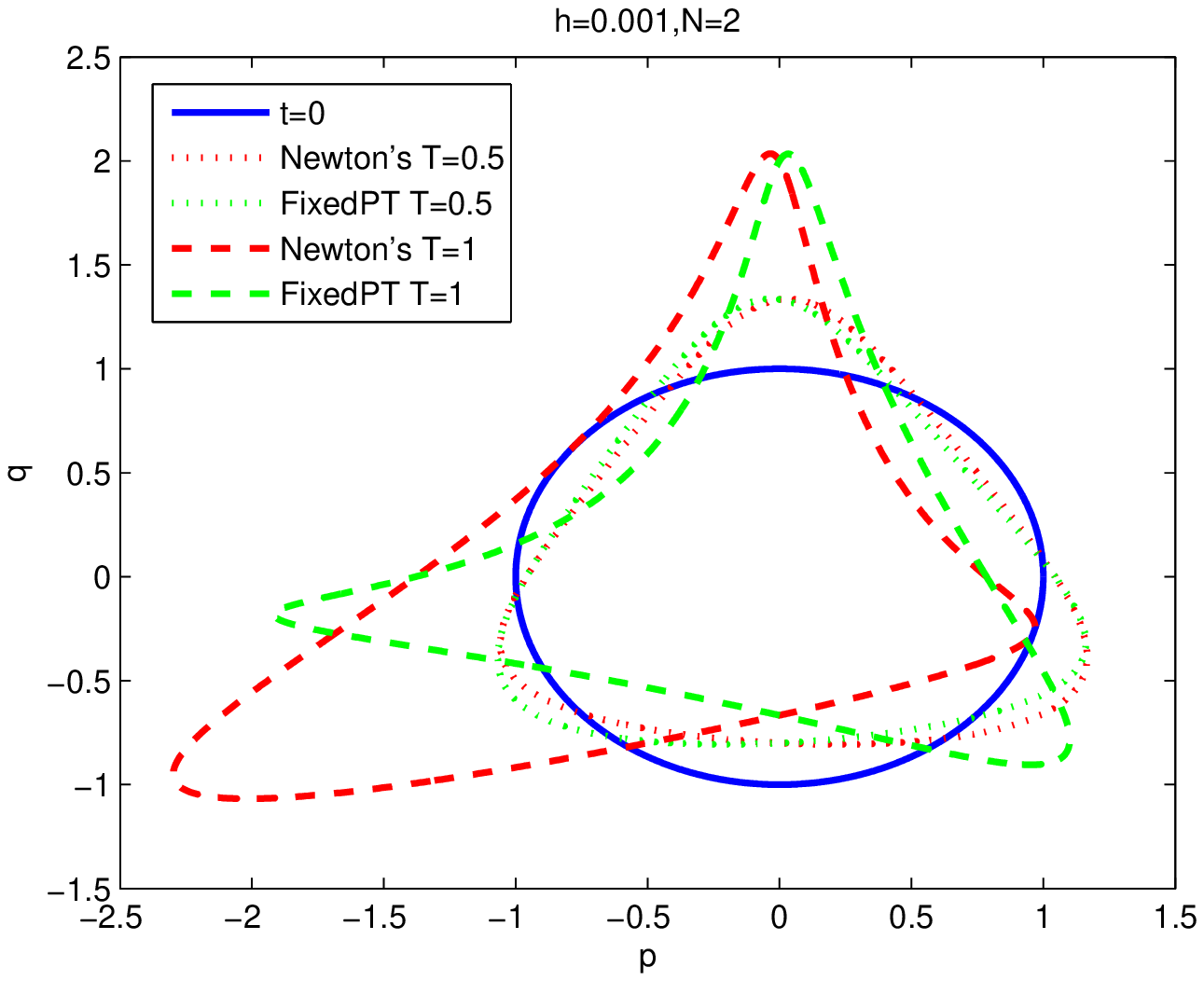}
\caption{Evolution of the unit circle driven by the midpoint rule (\ref{midshs}) with (left) Newton's iteration for $N=2$ (red dotted), $N=4$ (red dashed), fixed-point iteration for $N=2$ (green dotted), $N=4$ (green dashed), and $h=0.001$, $T=1$; (middle) Newton's iteration for $h=0.001$ (red dotted), $h=0.002$ (red dashed), fixed-point iteration for $h=0.001$ (green dotted), $h=0.002$ (green dashed), and $N=2$, $T=1$; (right) Newton's iteration for $T=0.5$ (red dotted), $T=1$ (red dashed), fixed-point iteration for $T=0.5$ (green dotted), $T=1$ (green dashed), and $N=2$, $h=0.001$.}\label{f9}
  \end{center}
\end{figure}

Figure \ref{f9} is devoted to the comparison between the Newton's iteration (\ref{newtshs}) and the fixed-point iteration
\begin{equation}\label{fixedshs}
\begin{split}
p_{n+1}^{[N]}&=p_n-h(\frac{p_n+p_{n+1}^{[N-1]}}{2})(\frac{q_n+q_{n+1}^{[N-1]}}{2})-\frac{1}{2}\overline{\Delta W_n}(\frac{p_n+p_{n+1}^{[N-1]}}{2})^2,\quad p_{n+1}^{[0]}=p_n,\\
q_{n+1}^{[N]}&=q_n+\frac{h}{2}(\frac{q_n+q_{n+1}^{[N-1]}}{2})^2+\overline{\Delta W_n}(\frac{p_n+p_{n+1}^{[N-1]}}{2})(\frac{q_n+q_{n+1}^{[N-1]}}{2}),\quad q_{n+1}^{[0]}=q_n,
\end{split}
\end{equation}
with $N=1,2,\cdots$.

Again, the left panel of Figure \ref{f9} indicates the little influence with mild change of $N$ for both iteration algorithms.
There are similar abilities of area preservation by the two iteration methods under the given data setting, while the middle panel of Figure \ref{f9} shows sensitivity of the Newton's iteration with respect to $h$, and the relative stability of the fixed-point iteration with respect to $h$.
As demonstrated by the right panel of Figure \ref{f9}, with the increase of the terminal time $T$, the error in area-preservation by both iteration methods grows, while that by the Newton's iteration seems more obvious, which maybe due to its more complicated calculations that are much easier to accumulate round-off errors.

\section{Concluding Remarks}
\label{7}
Preserving quadratic invariants is an advantage of numerical methods for performing good numerical behavior in convergence and stability.
Stochastic implicit Runge-Kutta methods can posses QI under certain conditions, while they are in general difficult to be realized directly.
But when we implement actually implicit SRK methods usually by taking place by either explicit approximates or iterations, which result in loss of accurate preservation of quadratic invariants. Based on the combinatory theory of rooted colored trees,
 this paper is devoted to  firstly give the conditions of preserving the quadratic invariants up to certain orders for explicit SRK methods.
Furthermore, this conditions enable establish the nearly conservative explicit SRK methods.
 Secondly we quantitatively analysis the loss for both explicit approximates and iterative implementation with fixed-point and Newton's iterations.
 The bounds of errors in the preservation of quadratic invariants by fixed-point and Newton's iterations are provided, which reveals the convergence rate of the errors with respect to the iteration number $N$, and the time step-size $h$.
  Numerical experiments are performed to testify the theoretical results, which support the theoretical results, and suggest appropriate choice of methods and parameters.
  Meanwhile, though with a theoretically faster convergence,
the Newton's iteration is shown in the numerical tests to be more critical and sensible for the choice of $h$ and the initial points than the fixed-point iteration.
In other words, the fixed-point iteration might be more reliable in the stochastic context.

\bibliographystyle{IMANUM-BIB}
\bibliography{IMANUM-refs}

\end{document}